\documentclass{article}
\usepackage{amsmath, amsthm, amssymb, setspace, graphics}
\usepackage[all]{xy}
\usepackage{url}
\newtheorem{thm}{Theorem}[subsection] 
\newtheorem{pro}[thm]{Proposition} 
\newtheorem{lem}[thm]{Lemma} 
\newtheorem{cor}[thm]{Corollary} 
\theoremstyle{definition} 
\newtheorem{defn}[thm]{Definition} 
\theoremstyle{remark} 
\newtheorem{rem}[thm]{Remark}
\newtheorem{notn}[thm]{Notation}
\newtheorem{exa}[thm]{Example}

\newcommand{\into}{\hookrightarrow}

\newcommand{\CC}{\mathbb C}

\newcommand{\PP}{\mathbb P}
\newcommand{\QQ}{\mathbb Q}

\newcommand{\ZZ}{\mathbb Z}

\newcommand{\Ocal}{\mathcal O}
\newcommand{\Ccal}{\mathcal C}
\newcommand{\Kcal}{\mathcal K}

\newcommand{\rest}[1]{{\textstyle|}_{#1}}

\newcommand{\Aut}{\operatorname{Aut}}

\newcommand{\Div}{\operatorname{Div}}
\newcommand{\pt}{\operatorname{pt}}

\newcommand{\Sing}{\operatorname{Sing}}

\newcommand{\hcf}{\operatorname{hcf}}

\newcommand{\Pic}{\operatorname{Pic}}

\newcommand{\Spec}{\operatorname{Spec}}

\newcommand{\Vect}{\operatorname{Vect}}
\newcommand{\al}{\alpha}
\newcommand{\be}{\beta}

\newcommand{\ep}{\varepsilon}

\newcommand{\la}{\lambda}

\newcommand{\ctop}{{\rm c}_{\mathrm{top}}}

\newcommand{\irr}{\operatorname{irr}}

\newcommand{\sh}{\operatorname{sh}}
\newcommand{\ch}{\operatorname{ch}}
\newcommand{\td}{\operatorname{td}}

\newcommand{\Sym}{{\operatorname{Sym}}}
\newcommand{\textsum}{{\textstyle{\sum}}}
\newcommand{\textsqcup}{\textstyle{\bigsqcup}}

\newcommand{\lvec}[1]{\vec{#1}}

\DeclareMathAlphabet{\mathpzc}{OT1}{pzc}{m}{it}
\newcommand{\sta}{\mathcal}
\def\pri#1{{#1}'}

\newcommand{\olC}{\overline C}

\def\pmmu{{\pmb \mu}}

\def\ol{\overline}

\def\wt{\widetilde}
\def\wh{\widehat}
\newcommand{\coa}{\ol}
\newcommand{\DDMM}{\overline{\sta M}}

\newcommand{\MMM}{\sta M}
\newcommand{\SSS}{\sta S}
\newcommand{\CCC}{{\sta C}}
\newcommand{\MMMbar}{\overline{\sta M}}
\begin{document}

\title{\textbf{Towards an enumerative geometry of the moduli
space of twisted curves and $r$th roots}}
\date{March 9, 2008}
\author{Alessandro Chiodo\footnote
{Partially supported by the
Marie Curie Intra-European Fellowship within the 6th European
Community Framework Programme, MEIF-CT-2003-501940.}} \maketitle

\begin{abstract}
The enumerative geometry of $r$th roots of line bundles is
crucial in the theory of $r$-spin curves and occurs in
the calculation of Gromov--Witten invariants of
orbifolds. It requires the definition of the suitable compact
moduli stack
and the generalization of the standard techniques from
the theory of moduli of stable curves.
In \cite{Ch_mod}, we constructed a compact stack by describing the
notion of stability in the context of twisted curves.
In this paper, by working with stable twisted curves,
we extend Mumford's formula for the Chern character of
the Hodge bundle to the direct image of the
universal $r$th root in $K$-theory.
\end{abstract}

\maketitle

\vspace*{6pt}\tableofcontents  

\section{Introduction}
\label{sec:introduction}
In this paper we investigate the enumerative
geometry of the moduli of curves paired with the $r$th roots of a
given line bundle (usually the canonical bundle or the structure sheaf).
It is well known that
one can construct a proper and smooth stack
$\MMMbar^{r}_{g,n}$ of $r$th roots of the canonical bundle
of genus-$g$ curves. More generally, one can consider $r$th roots of
any power of
the canonical bundle possibly twisted at the markings $x_1, x_2, \dots, x_n$ with
arbitrarily chosen  multiplicities $m_1, m_2, \dots, m_n$; see \cite{Ja_geom}, \cite{AJ},
\cite{CCC}, and \cite{Ch_mod}.
In Theorem \ref{thm:GRRcalc},
we determine the Chern character
$$\ch(R^\bullet\pi_*\SSS),$$
where $\pi\colon \CCC\to \MMMbar_{g,n}^{r}$ is the universal curve and
$\SSS$ is the universal $r$th root.
As illustrated in Section \ref{sect:appl},
the formula for $\ch(R^\bullet\pi_*\SSS)$
has direct applications to
Gromov--Witten theory of stacks,
to $r$-spin theory in the sense of
Witten \cite{Wi} (see also Fan, Jarvis, and Ruan \cite{FJR}
for recent generalizations and conjectures), and
to the mathematics of
Hurwitz numbers (see Zvonkine \cite{ZvLu}).

\subsection{Main theorem: Grothendieck Riemann--Roch for universal $r$th roots}
\label{sect:intromainthm}
Via Grothendieck Riemann--Roch, in Theorem \ref{thm:GRRcalc}, we express the
Chern character of the perfect complex
$R^\bullet \pi_* \SSS$
in terms of tautological classes: the kappa classes $\kappa^{}_d$
and the psi classes $\psi_1, \dots, \psi_n$.
The calculation takes place on
the stack $\MMMbar^{r}_{g,n}$, the proper moduli space of
$r$th roots of the line bundle $\mathcal K=(\omega^{\log})^{\otimes s}(-\textsum_{i=1}^n m_i[x_i])$,
where $s$, $m_1, \dots, m_n$ are integers, and we adopt the notation
$$\omega^{\log}=\omega(\textsum_{i=1}^n [x_i]).$$ The kappa and psi classes
$$\kappa^{}_d=\pi_* ({\rm c}_1(\omega^{\log})^{d+1}) \quad
\quad \text{and} \quad \quad \psi_i={\rm c}_1(\omega\rest{\text{marking $x_i$}}),$$
are defined (see \eqref{eq:defkappapsi}) in the Chow ring of $\MMMbar^{r}_{g,n}$ rather than
in the Chow ring of
$\MMMbar_{g,n}$ as it usually occurs.
In view of applications to enumerative geometry, we
provide a formula which  allows us to compute intersection numbers
directly
via products in the tautological ring of $\MMMbar_{g,n}$,
Corollary \ref{cor:pushdown}.

\paragraph{The case $r=1$: Mumford's calculation for the Hodge bundle.}
The Hodge bundle $\mathbb H$ is the direct image  of $\omega$ on $\MMMbar_{g,n}$.
In \cite{Mu_GRR},
Mumford provides a formula for $\ch(\mathbb H)$ in terms of tautological classes
and Bernoulli numbers.
In view of our generalization, it is convenient to reformulate Mumford's formula in terms of
Bernoulli polynomials $B_m(x)$ by recalling the relation
$B_m(1) =(-1)^m B_m,$ where $B_m$ is the usual Bernoulli number.
We rewrite Mumford's formula as follows
\begin{equation}\label{eq:Mu_GRR}
\ch(R^\bullet \pi_* \omega)=\ch(\mathbb H)-1=\sum_{d\ge 0} \left(\frac{B_{d+1}(1)}{(d+1)!}\kappa^{}_d -\sum _{i=1}^n
\frac{B_{d+1}(1)}{(d+1)!}\psi_i^d+\frac{1}{2}\frac{B_{d+1}(1)}{(d+1)!}j_*(\gamma_{d-1})\right),
\end{equation}
where $j$ is the composite of $\pi$ and
the natural double cover
$$\pri{\Sing}\longrightarrow \Sing\subset {\CCC}$$
(we recall that the singular locus $\Sing$
of ${\CCC}\to \MMMbar_{g,n}$
consists of nodes
and the category $\Sing'$ formed by nodes
equipped with the choice of an order for the
two branches of each node is a double cover).
The class $\gamma_{d}\in A^{d}(\Sing')$ is
defined as follows. Let $\psi$ and $\widehat \psi$ be the first
Chern classes of the line bundles whose fibre at each point $p'$ of $\Sing'$
is the cotangent line along the first and the second branch, respectively. Then,
we set
$$\gamma_{d}:=\textsum_{i+j=d} (-\psi)^{i} \widehat \psi^j,$$
and we let $\gamma_d$ be zero for $d<0$.

Mumford's equation \cite[(5.2)]{Mu_GRR} for $\ch(\mathbb H)$
is written in a different form, but the above
expression can be obtained
once the following remarks are made:
\begin{enumerate}
\item Mumford does not use the standard notation for Bernoulli numbers recalled at \ref{notn:bernoulli}.
\item The definition of the kappa classes
recalled above is the standard notation for marked curves and differs from
the one of \cite{Mu_GRR}, which does not involve the markings: in \cite{Mu_GRR}
the kappa class is $\pi_* ({\rm c}_1(\omega)^{d+1})$ rather than
$\pi_* ({\rm c}_1(\omega^{\log})^{d+1})$. For example
Mumford's $0$th kappa class is $2g-2$ rather than
$\kappa_0=2g-2+n$. In degree $d$ this change makes
the term $-\sum_i \psi_i^d$ appear in our formula (see \eqref{eq:kappalog}).
\item Mumford has a typo in the index of the kappa term of his formula (the index ``$2l+1$''
should rather be ``$2l-1$'').
\end{enumerate}

\paragraph{The case $r>1$: stack-theoretic curves and Bernoulli polynomials.}
Our generalization of Mumford's theorem
takes place on the proper moduli stack $\MMMbar^{r}_{g,n}$
of $r$th roots of $$\mathcal K=(\omega^{\log})^{\otimes s}(-\textsum_{i=1}^n m_i [x_i]),$$
for $s, m_1,\dots, m_n\in \ZZ$ and $(2g-2+n)s-\sum_i m_i\in r\ZZ$. The formula
involves Bernoulli polynomials evaluated at $s/r$ and $m_i/r$ (see Notation \ref{notn:bernoulli}). In this way,
Mumford's calculation corresponds to the case where
$r, s,$ and the indices $m_i$ are all equal to $1$ and
$\MMMbar^{r}_{g,n}$ equals $\MMMbar_{g,n}$.

In \cite{Ch_mod}, we show  that the construction of $\MMMbar^{r}_{g,n}$
is just a slight modification of the
compactification $\MMMbar_{g,n}$ of
$\MMM_{g,n}$ which was obtained by Deligne and Mumford by allowing \emph{stable curves}
(curves $C$ with nodes as singularities
and finite automorphism group $\Aut(C)$).
The stack $\MMMbar^{r}_{g,n}$ is also a compactification: it
contains the stack $\MMM_{g,n}^r$ classifying
smooth curves with an $r$th root of $\mathcal K$.
It is also obtained by enlarging
the category $\MMM_{g,n}^r$
to new geometric objects:
$r$th roots of $\mathcal K$ on \emph{$r$-stable curves}, which are
stack-theoretic curves $\mathcal C$ with stable coarse space, representable smooth locus,
and nodes $p\colon \Spec \CC\to \mathcal C$ such that $\Aut_{\Ccal}(p)$ is
isomorphic to $\pmmu_r$.

As in Mumford's calculation, we need to consider
$\Sing$, the stack classifying curves (equipped with an $r$th root)
alongside with the choice of a singularity.
As above, the stack $\Sing$ naturally maps to $\MMMbar^r_{g,n}$ and has
a natural double cover induced by the choice of a  branch for each node;
in this way we have the morphism
$$j\colon \pri{\Sing}\to \MMMbar^{r}_{g,n}.$$
The stack $\pri{\Sing}$
is naturally equipped with two line bundles
whose fibres are the cotangent lines to the first branch of the
coarse curve and the cotangent lines to the second branch.
We write $\psi,\wh\psi \in H^2( \pri{\Sing}, \QQ)$
for their respective
first Chern classes.
(In this notation, we privilege the
coarse curve, because in this way the classes $\psi$
and $\wh\psi$ are more easily related to the
classes $\psi_i$.)

Recall that $\CCC$ is equipped with the universal $r$th root $\SSS$.
For each point $p'$
of ${\Sing}'$, the restriction of $\SSS$ to the
first branch yields a multiplicity index $q(\pri p)\in \{0, \dots, r-1\}$
as follows.
The line bundle
$\SSS$ on the point lifting the node to the first branch
can be regarded as a representation of $\pmmu_r$.
This representation is a power of the
representation attached to the line
cotangent to the first branch. We define $q(p')$ as the order of that power.
By construction, such a multiplicity index is locally constant on $\pri{\Sing}$.
We get a decomposition of $\pri{\Sing}$ and the restriction morphisms
$$\pri{\Sing}={\textstyle \bigsqcup _{q=0}^{r-1}\pri{\Sing}_q}\quad \quad\quad \quad  \quad \quad
j_q=j\rest {\pri{\Sing}_q}\colon {\Sing}'_q\to \MMMbar^{r}_{g,n} \quad \text{for }0\le q<r.$$

\begin{thm}\label{thm:GRRcalc}
Let $s, m_1, \dots, m_n$ be integers satisfying $(2g-2+n)s-\sum_i m_i\in r\ZZ$.
Let $\SSS$ be the universal $r$th root of
$(\omega^{\log})^{\otimes s}(-\textsum_{i=1}^n m_i[x_i])$
on the universal $r$-stable curve $\CCC$ over the moduli stack
$\MMMbar_{g,n}^r$.
The direct image $R^\bullet\pi_*\SSS$ via the universal curve
$\pi\colon \CCC\to \MMMbar_{g,n}^r$ satisfies the
equation
$$\ch(R^\bullet \pi_* \SSS)=
\sum_{d\ge0}\left(\frac{B_{d+1}\left(s/r\right)}{(d+1)!} \kappa^{}_d -
\sum_{i=1}^n \frac{B_{d+1}\left(m_i/r\right)}{(d+1)!}\psi_i^{d}+
\frac{1}{2}\sum_{q=0}^{r-1}\frac{rB_{d+1}(q/r)}{(d+1)!}(j_q)_*
\left({\gamma_{d-1}}\right)\right),$$
where the cycle $\gamma_{d}$
in $A^{d}(\pri \Sing_q)$ is
$\textsum_{i+j=d} (-\psi)^{i} \widehat \psi^j$ (we set $\gamma_d=0$ when $d$ is negative).
\end{thm}

\subsection{Applications}\label{sect:intromotiv}
Theorem \ref{thm:GRRcalc} can be used to compute several
invariants in the enumerative geometry of curves. In Section \ref{sect:appl}, we
illustrate the relation with the crepant resolution conjecture
in Gromov--Witten theory and with the theory of $r$-spin curves.

\paragraph{The crepant resolution conjecture.}
The first context where Theorem \ref{thm:GRRcalc} provides
a useful application is
the Gromov--Witten theory of orbifolds.
Indeed, when we set $s=0$, Theorem \ref{thm:GRRcalc} can be regarded as a tool for the
computation of genus-$g$ Gromov--Witten invariants
of the orbifold $[\CC^2/\pmmu_r]$,
where $\pmmu_r$ acts as $(x,y)\mapsto(\xi_r x,\xi_r^{-1}y)$.
In section \S\ref{sect:applCRC}, we detail the definition of such invariants,
we spell out how Theorem \ref{thm:GRRcalc} applies, and
we provide some examples, which
can be generalized for any value of $g$ and $r$.
The crucial observation is that stable maps to $[\CC^2/\pmmu_r]$
are related to $r$th root of $\Ocal(-D)$, where $D$ is a divisor supported on the markings.

The Gromov--Witten potential is the power series assembling all Gromov--Witten invariants.
The potential of
Gorenstein orbifolds such as $[\CC^2/\pmmu_r]$ is the subject of
the crepant resolution conjecture.
In the statement of Bryan and Graber \cite{BG}, the conjecture
predicts that the genus-zero potential
of a Gorenstein orbifold and that of any
crepant resolution of the coarse space are equal after
a change of variables.
Recently, this statement has been proved for the genus-zero case by Corti, Coates, Iritani and Tseng \cite{CCIT}
and generalized to a conjecture involving higher genera by Coates and Ruan \cite[Conj. 10.2]{CR}.

\paragraph{Invariants in $r$-spin theory.}
The theory of $r$-spin curves is the next natural context where
Theorem \ref{thm:GRRcalc} applies. By definition, $r$-spin structures
are $r$th roots of $\omega(-D)$, where $D$ is
an effective divisor supported on the markings. This corresponds to the
case $s=1$ (recent work by \cite{FJR} cast this notion  in the
wider  framework of $W$-spin theory; we remark that this
new notion directly involves higher values of $s$).

The crucial invariant of this theory is again a power series, the $r$-spin potential,
assembling several intersection numbers in $\MMMbar_{g,n}^r$ (see discussion in \S\ref{sect:applWspin}).
This theory has been put in the framework of Gromov--Witten
theory by Jarvis, Kimura, and Vaintrob \cite{JKV2}.
In \cite{FSZ}, Faber, Shadrin, and Zvonkine prove that the $r$-spin potential satisfies the
Gelfand--Diki\u\i\ hierarchy as conjectured by Witten in \cite{Wi}.
In \cite{CZ}, in collaboration with Zvonkine,
Theorem \ref{thm:GRRcalc} is used to provide a formula relating this
Gromov--Witten $r$-spin potential to the
twisted Gromov--Witten $r$-spin potential (the $r$-spin potentials
extended with the Chern characters $\ch(R^\bullet\pi_*\SSS)$).
Finally, Theorem \ref{thm:GRRcalc} provides new information on the enumerative geometry of
$R^\bullet \pi_* \sta S$; we discuss in some examples the relation with the
cycle ${\rm c}_{\rm W}$, the virtual fundamental cycle in Gromov--Witten
$r$-spin theory.

We also point out that the Chern classes of
$R^\bullet \pi_* \sta S$ are expected to be related
to Hurwitz numbers via a conjecture by Zvonkine (see \cite{ZvLu}),
which generalizes the ELSV formula \eqref{eq:ELSV} concerning the
Chern classes of the Hodge bundle.

\paragraph{GRR formulae.} The proof of
Theorem \ref{thm:GRRcalc} provides a new application of
the Grothendieck Riemann--Roch formula.
In Section \ref{sect:GRR},
we start from a stack-theoretic curve $\wt C$ equipped with
the $r$th root $\wt S$ on a base scheme $X$; then, we work in the cohomology ring of the
desingularization $C$
of the coarse space $\ol C$ of $\wt C$. In this way,
the curve $C$ is  semistable and is
equipped with a line bundle $S$ whose direct image on the base scheme $X$ coincides with
that of the universal $r$th root $\wt S$ on $\wt C$, Lemma \ref{lem:Dq}.
This allows us to carry out the calculation in the scheme-theoretic context via the
standard GRR formula for schemes.

The statement of Theorem \ref{thm:GRRcalc} can also be obtained
via Toen's GRR
 formula for  stacks \cite{To}:
this alternative approach is illustrated in \S\ref{sect:staGRR}
and the result coincides with the formula of Theorem \ref{thm:GRRcalc}.
In this comparison, the main difficulty
is the identification of the formula of Theorem \ref{thm:GRRcalc}, which lies in rational cohomology,
with the formula given by Toen's theorem, which is naturally stated
in terms of cohomology classes with ``coefficients in the representations".
We point out that
this problem also
occurs  in the proof of
Tseng's orbifold quantum Riemann--Roch, the
formula expressing
the generating function of twisted Gromov--Witten invariants
in terms of the generating function
of the untwisted invariants, \cite{Ts}. There, stack-theoretic curves
also occur, and Toen's GRR is applied; in \cite[\S7.2.6]{Ts},
by means of calculations in the $r$th cyclotomic field,
the ``coefficients in the representations''
are expressed in terms of Bernoulli polynomials.

\subsection{Terminology}
Throughout this paper we work over the field of complex numbers $\CC$.

\begin{notn}
By $r$ we always denote a positive integer. We often
use the involution of $\{0, \dots, r-1\}$ sending $q>0$ to $r-q$ and
fixing $0$. We write it as $q\mapsto \ol{r-q}$,
by  adopting
the convention $\ol i = i$ if $0\le i<r$ and $\ol r=0$.
\end{notn}

\begin{notn}
We adopt the terminology of
\cite{LM} for stacks.
For Deligne--Mumford stacks there is a natural $2$-functor to the category of
algebraic spaces associating
to a stack
an algebraic space which is universal with respect to morphisms to algebraic spaces.
We usually refer to it as the {coarse space}.
\end{notn}

\begin{notn}
We often need to describe stacks and morphisms between stacks locally
in terms of explicit equations.
Let ${X}$ and $U$ be algebraic stacks and
let $x\in X$ and $u\in U$ be geometric points.
We say ``the local picture of ${X}$
at $x$ is
given by ${U}$ (at $u$)'' if
there is an isomorphism between the
strict henselization ${X}^{\text{sh}}$
of ${X}$ at $x$ and the strict henselization
${U}^{\sh}$
of ${U}$ at $u$.
\end{notn}

\begin{notn} Let $f\colon X\to Y$ be a proper birational morphism of complex varieties.
We write ${\rm Ex}(f)$ for the exceptional locus in $X$ of the
morphism $f$.
\end{notn}

\begin{notn}[moduli of stable curves]
For $2g-2+n>0$, the category $\DDMM_{g,n}$ of stable curves of genus $g$ with $n$ markings
is a proper, smooth, and irreducible stack of Deligne--Mumford type \cite{DM} of dimension
$3g-3+n$.
The full subcategory $\MMM_{g,n}$ of smooth curves
is an open substack.
\end{notn}

\begin{rem}\label{rem:omega}
We often consider nodal curves $C\to X$ with sections $s\colon X\to C$ in the smooth locus.
We point out that there is a canonical isomorphism between
$s^*\omega_{C/X}(s(X))$ and $\Ocal _X$.
\end{rem}

\subsection{Structure of the paper}
In Section \ref{sect:setup} we recall the construction of $\MMMbar_{g,n}^r$ and
we set up the notation.
In Section \ref{sect:GRR} we prove the theorem stated above.
In Section \ref{sect:appl},
we illustrate the applications mentioned in \eqref{sect:intromotiv}.

\subsection{Acknowledgements}
I would like to thank D.~Zvonkine, for pointing out
Example \ref{exa:torsion1} to me and
for showing me his conjectural statement relating $r$-spin
curves and Hurwitz numbers.

In the course of the proof of Theorem \ref{thm:GRRcalc}, we
handle the degeneration of roots  by means of $q$-chain divisors. These
divisors are implicitly used by L.~Caporaso, C.~Casagrande, and M.~Cornalba
in \cite{CCC}, and I am grateful to M.~Cornalba for illustrating
them to me (see Fig.~3).

I would like to thank J.~Bertin and the anonymous referee for the care
with which they read this work and for
their precious comments and corrections.
My thanks to S.~Boissi\`ere, N.~Borne, T.~Graber, F.~Herbaut,
A.~Hirschowitz, T.~Rivoal, and C.~Simpson
for helpful conversations and advice.

\section{Previous results on moduli of curves and $r$th roots}\label{sect:setup}
In this section, we recall the definition of $\MMMbar_{g,n}^{r}$, the
proper moduli stack of $r$th roots of
$$\mathcal K=(\omega^{\log})^{\otimes s}(\textsum_{i=1}^n (-m_i[x_i]),$$
for $s, m_1,\dots, m_n\in \ZZ$.
For degree reasons, the category of $r$th roots of $\mathcal K$ is
nonempty if and only if we have
$$s(2g-2+n)-\textsum_{i=1}^n m_i\in r\ZZ.$$
It is known that there exists a proper stack of $r$th roots of $\mathcal K$,
see \cite{Ja_geom},
\cite{AJ}, \cite{CCC}, and \cite{Ch_mod}. Here we follow the construction
of \cite{Ch_mod}, where we show that the only difficulty lies in
defining a new compactification of $\MMM_{g,n}$.
\subsection{Compactifying $\MMM_{g,n}$ using $r$-stable curves
instead of stable curves}
\begin{defn}[$r$-stable curve]
\label{defn:rstable}
For $X$ a scheme and $C$ a stack, consider
a proper and flat morphism
$C\to X$ with $n$ ordered section $x_1,\dots, x_n \colon X\to C$,
mapping into the smooth locus of $C\to X$
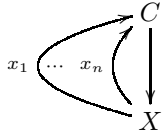
\[ \xymatrix@R=1cm{
{C}\ar[d]\\
X\ar@/^1.5cm/[u]^{x_1}_\cdots\ar@/^0.5cm/[u]^{x_n}\\
} \]
The data $(C\to X, x_1,\dots,x_n)$ form
an  \emph{$r$-stable} curve
of genus $g$ with $n$ markings
if
\begin{enumerate}
\item
the fibres are purely $1$-dimensional
with at most nodal singularities,
\item the smooth locus
$C^{\rm sm}$ is an algebraic space,
\item the coarse space $\coa{C}\to X$ with the
sections $\coa {x}_1,\dots, \coa{x}_n$
is a genus-$g$  $n$-pointed stable curve
$(\coa{C}\to X, \coa{x}_1,\dots,\coa{x}_n)$,

\item
the local picture
at a node
is given by $[U/\pmmu_r]\to T$, where
\begin{enumerate}
\item[$\bullet$] $T=\Spec A$,
\item[$\bullet$] $U=\Spec A[z,w]/(zw-t)$ for some $t\in A$, and
\item[$\bullet$] the action of $\pmmu_r$ is given by
$(z,w)\mapsto (\xi_rz,\xi_r^{-1}w)$.
\end{enumerate}
\end{enumerate}
\end{defn}

\begin{rem}
For any $r$-stable curve as above
we denote by $\omega_{C/X}^{\log}$
the sheaf of logarithmic differentials
defined as
\begin{equation}\label{eq:omegalog}\omega_{C/X}^{\log}:=
\omega_{C/X}\otimes \Ocal(\textsum_{i=1}^n x_i(X)).
\end{equation}
If
$(\coa{C}\to X, \coa{x}_1,\dots,\coa{x}_n)$ is the corresponding
coarse stable curve, we have \cite[Prop.~2.5.1]{Ch_mod}
\begin{equation}\label{eq:omegapull}
\pi_C^*(\omega^{\log}_{\coa{C}/X})\cong\omega^{\log} _{C/X}.
\end{equation}
\end{rem}

For each $r$, the notion of $r$-stable
curve yields a compactification of the
stack of $n$-pointed genus-$g$ stable smooth curves
(recall that stability in the case
of smooth curves only means $2g-2+n>0$).
The following theorem is a consequence of a theorem of Olsson \cite{Ol}
and is stated in this form in {\cite[4.1.4, 4.2.7]{Ch_mod}.
\begin{thm}[\cite{Ol}, \cite{Ch_mod}]
The moduli functor $\DDMM_{g,n}(r)$ of
$r$-stable $n$-pointed curves of genus $g$ forms a proper, smooth, and
irreducible Deligne--Mumford stack of dimension $3g-3+n$.

The stack of $1$-stable curves
$\DDMM_{g,n}(1)$ is the stack of stable curves ${\DDMM}_{g,n}$ in
the sense of Deligne and Mumford.
For any $r> 1$ and for $2g-2+n>0$,
we get new compactifications of the stack of smooth $n$-pointed curves
$\MMM_{g,n}$. There is a natural
surjective, finite, and flat morphism
$\DDMM_{g,n}(r)\to \DDMM_{g,n}(1)$ which is invertible
on the open dense substack of smooth $n$-pointed genus-$g$ curves
and yields an isomorphism between coarse spaces.
The morphism is not injective; indeed, the restriction of $\DDMM_{g,n}(r)\to \MMMbar_{g,n}$
to the substack of singular $r$-stable curves
has degree $1/r$ as a morphism between stacks.
\qed
\end{thm}

\begin{rem}[$\vec{l}$-stability]\label{rem:lstab}
In fact in \cite{Ch_mod}, we classify all compactifications of $\MMM_g$ inside
Olsson's nonseparated stack of twisted curves: stack-theoretic curves with
representable smooth locus
and nodal singularities with finite stabilizers. In this way, we hit upon the
notion of \emph{$\lvec{l}$-stability}, which we recall in passing even if it will not be used
in this paper. For simplicity, let us consider the case of unmarked curves.
 A node of a stable curve over $\CC$
has type $0$ if the normalization of the stable curve at
the node
is connected, whereas it has type
$i\in \{1, \dots,{\lfloor g/2\rfloor }\} $ if the normalization
is the union of two disjoint curves
of genus $i$ and $g-i$. For $\vec{l}=(l_i\mid i\in \{0, 1,\dots,{\lfloor g/2\rfloor }\})$,
a twisted curve is $\vec{l}$-stable if its nodes of type $i$ have stabilizers of order $l_i$.
We point out that when all entries of $\vec{l}$ equal $r$, an $\vec{l}$-stable curve
is nothing more than an $r$-stable curve in the sense of Definition \ref{defn:rstable}.
\end{rem}

\begin{rem}[the singular locus of the universal $r$-stable curve]
\label{rem:singloc}
In Section \ref{sect:GRR}, we work in the \'etale topology of $\DDMM_{g,n}(r)$.
Let $X$ be a scheme and $X\to \MMMbar_{g,n}(r)$
be an \'etale morphism; then, the corresponding $r$-stable
$X$-curve is a morphism $\wt \pi\colon \wt C\to X$ whose variation is maximal:
the Kodaira--Spencer homomorphism
$$T\MMMbar_{g,n}(r)_{x}\to  \text{Ext}^1(\sta C_x,\Ocal_{\sta C_x})$$
is surjective for any geometric point $x$ in $\MMMbar_{g,n}(r)$.
In this case, the following facts derive from the theory of
moduli of nodal curves and $r$-stable curves (see \cite{Ol} and \cite{Ch_mod}).
\begin{enumerate}
\item[A)] {The stack $\wt C$:} \newline
The $r$-stable curve $\wt C$ is a regular stack and the
singular locus $\Sing(\wt C/X)$
of $\wt \pi\colon \wt C\to X$
is regular  of codimension two in
$\wt C$. The image of $\Sing(\wt C/X)$
in $X$ with
respect to $\pi$ is a normal crossing divisor (the pullback in  $X$ of
the boundary locus  via $X\to \MMMbar_{g,n}(r)$, see \cite{Ol}).
Note that the local picture of $\wt \pi$ at
any point $x\in \Sing(\wt C/X)$ is
given by
$$[\Spec U/\pmmu_r] \to \Spec T,$$
where $T$ equals $\CC[x_1,\dots, x_{3g-3+n}]$, $U$ equals $T[z,w]/(zw=x_1)$, and the group $\pmmu_r$ acts 
as
$(z,w)\mapsto (\xi_r z, \xi_r^{-1} w)$.
\item[B)]{The coarse space $\olC$ of $\wt C$:} \newline
If we consider the coarse space of $\wt C$; then, we get
an algebraic scheme $\overline{C}$ with singularities
of  type $A_{r-1}$ on every point of the singular locus.
More precisely, using the local picture given at the previous point,
we see that the singular locus $\Sing(\olC/X)$ of the stable curve
$\overline{\pi}\colon \olC \to X$ is a regular scheme, has codimension two in $\olC$, and
the local picture of $\overline{\pi}$ at any of its points is
given by the spectrum of
$$\CC[x_1,\dots, x_{3g-3+n}]\longrightarrow
\CC[x_1,\dots,x_{3g-3+n}, z, w]/(zw=x_1^r).$$
Note that $\wt\pi(\Sing(\wt C/X))$ equals $\overline{\pi} (\Sing(\olC/X))$
and is a normal crossings divisor in $X$.
\item[C)] {The desingularization $C$ of $\olC$:} \newline
\begin{figure}[hh]
\begin{picture}(200,70)(110,-10)
  \qbezier[30](395,5)(405,40)(445,50)
  \qbezier[30](185,5)(175,40)(135,50)
  \qbezier(170,5)(215,50)(215,50)
  \qbezier(200,50)(245,5)(245,5)
  \qbezier(230,5)(275,50)(275,50)
  \qbezier(260,50)(305,5)(305,5)
  \put(295,30){$...$}
  \qbezier(410,5)(365,50)(365,50)
  \qbezier(380,50)(335,5)(335,5)
  \qbezier(350,5)(305,50)(305,50)
  \put(450,30){\vector(1,0){30}}
  \put(485,30){$\pt\in \Sing(\ol C/X)$}
 \end{picture}
\caption{The fibre of ${\rm Ex}(C\to \ol C)\to \Sing(\ol C/X)$, a chain of $r-1$ rational lines.}
\end{figure}

By repeating $\lfloor r/2\rfloor$ blowups of the singular locus
starting from $\overline{C}$ we get the desingularization $\pi\colon C\to \ol{C}$
of $\overline{\pi} \colon \olC\to X$. The morphism $\pi\colon C\to X$  is a semistable curve. The
fibre of $\pi$ is sketched in the figure (the exceptional locus
of $C\to \ol C$ is drawn with continuous lines and the
rest of the curve with dotted lines).
The scheme $C$ is regular, and the singular locus
$\Sing(C/X)$ is a regular scheme, of codimension two inside $C$,
and in fact it is an \'etale cover of $\Sing(\olC/X)$ of degree $r$.
Note that, for $r>1$, $\Sing(C/X)$ lies in the exceptional locus
${\rm Ex}(C\to \ol C)$ of
the desingularization $\nu \colon C\to \olC$ (this happens because for $r> 1$ every node of
the curve $\olC\to X$ is a singular point of $\olC$).

\end{enumerate}
\end{rem}

We now show how the compactification $\MMMbar_{g,n}(r)$ illustrated above
allows us to define easily a proper stack of $r$th roots of $\mathcal K$.
The above description of the singular loci will be generalized and
further developed in Remark \ref{notn:charclsing}.

\subsection{Roots of order $r$ on $r$-stable curves}
Once the moduli stack of $r$-stable curves is defined,
the construction follows naturally.
We fix the parameters $s, m_1, \dots, m_n\in \ZZ$ and we assume that $(2g-2+n)s-\sum _i m_i$ is a multiple of $r$.
The category of $r$th roots of $$\mathcal K=(\omega^{\log})^{\otimes s}(-\textsum_i m_i [x_i]),$$
on \emph{stable} curves with $n$ markings $x_i\in C$ forms a stack which is not proper;
instead, we consider $r$th
roots of $\mathcal K$ on \emph{$r$-stable} curves.
\begin{thm}[{\cite[Thm.~4.2.3]{Ch_mod}}]\label{thm:etale}
The moduli functor $\MMMbar^{r}_{g,n}$
of $r$th roots of
$\mathcal K$ on $r$-stable curves forms a proper and smooth stack of Deligne--Mumford type
of dimension $3g-3+n$.

For $s=m_1=\dots=m_n=0$, $\mathcal K$ is $\Ocal$ and the stack
$\MMMbar_{g,n}^{r}$
is a group stack $\sta G$
on $\MMMbar_{g,n}(r)$.
In general  $\MMMbar_{g,n}^{r}$ is a finite torsor on
$\MMMbar_{g,n}(r)$ under $\sta G$.

The morphism $\MMMbar_{g,n}^{r}\to \MMMbar_{g,n}(r)$
is \'etale. It factors through a morphism
locally isomorphic to $B\pmmu_r\to \Spec \CC$ (a
$\pmmu_r$-gerbe) and a representable
\'etale $r^{2g}$-fold cover; therefore
it is has degree $r^{2g-1}$.\qed
\end{thm}
\begin{rem}
In fact, in \cite[Thm.~4.2.3]{Ch_mod} we completely
determine the values of $\lvec{l}$ for which we
have the properness of the stack
of $r$th roots of $\mathcal K$
on curves satisfying $\vec{l}$-stability in the sense
illustrated in Remark \ref{rem:lstab}. As explained above, when all
entries of the multiindex $\vec{l}$ equal $r$, saying that a curve is
$r$-stable is equivalent to saying that a curve is $\vec{l}$-stable.
From the point of view of enumerative geometry, all choices of $\vec{l}$
are equivalent, because the rational cohomology of a Deligne--Mumford stack only
depends on the coarse space, and the coarse space of the compactification
does not depend on $\vec{l}$.
\end{rem}

\begin{rem}[analysis of the singular locus]\label{notn:charclsing}
As we did for the stack of $r$-stable curves
$\MMMbar_{g,n}(r)$ we study an object of $\MMMbar_{g,n}^r$
over a base scheme $X$ mapping to $\MMMbar_{g,n}^r$ via an \'etale
morphism. Via pullback we get an $r$-stable curve
$\wt\pi\colon \wt C\to X$ and an
$r$th root $\wt S$ of
$(\omega_{C/X}^{\log})^{\otimes s}(-\sum_i m_i [x_i])$.
These data involve several local geometric structures,
for which we  fix a notation for the rest of the paper.
\begin{enumerate}
\item[a)]
We write $$\wt Z=\Sing(\wt C/X)\into \wt C$$ for
the singular locus of the $r$-stable curve
$\wt\pi\colon \wt C\to X$.
Let $$\wt \varepsilon \colon
\wt Z'\to\wt Z$$ be the $2$-fold \'etale cover classifying the two possible
ways of ordering the incident branches.
We denote by $$\wt \iota\colon \wt Z'\to \wt Z',$$
the natural involution switching the order.
In this way, if we pullback $\wt C\to X$ via $\wt j\colon \wt Z'\to X$,
\[ \xymatrix@R=.1cm{
{\wt C\times_X \wt Z'}\ar[dd]\ar[rr]&&\wt C\ar[dd]\\
&\square &\\
\wt Z' \ar[rr] && X \\
} \]
we get a curve $\wt C\times_X \wt Z'\longrightarrow \wt Z'$ such that over each point
of $\wt Z'$ there is a prescribed choice of a node and of a branch of that node.

We use this fact in order to show that the stack $Z'$
admits a natural decomposition $$\wt {Z}'=\textsqcup_{q=0}^{r-1}
\wt{Z}_q'.$$
To this effect, we define a multiplicity index $q\colon \wt Z'\to \{0,\dots, r-1\}$
and we realize $\wt Z_0', \dots, \wt Z_{r-1}'$ as the loci where this index is constant.
Consider
a geometric point of $p'\in \wt Z'$ and the corresponding node
in the $r$-stable curve $\wt C\times_X \wt Z'\longrightarrow \wt Z'$. The
line cotangent to the prescribed branch is a line bundle on the node: a
faithful $\pmmu_r$-representation.
The restriction of the $r$th root $\wt S$
to the node can also be regarded as a $\pmmu_r$-representation: a power
of order $q\in \{0,\dots, r-1\}$ of the cotangent line to the prescribed branch.
In this way, the locally constant function $q$ is defined, and so are
 $\wt{Z}_0',\dots, \wt Z'_{r-1}$.
\item[b)]
Let $\ol \pi\colon \ol C\to X$ be the coarse space of $\wt \pi\colon \wt C\to X$.
We write $$Z=\Sing(\ol C/X)\into \ol C$$ for
the singular locus of the stable curve $\ol C\to X$.
As above, set $$\ol \ep \colon Z'\to Z \qquad
\text{and} \qquad \ol \iota\colon Z'\to Z'$$
 for the
natural $2$-fold \'etale cover of $Z$ and the involution of $Z'$.
By $$\ol j\colon Z'\to X,$$ we denote
the composite of $Z\to \ol C$ and $\ol C\to X$.
We write
$$\tau ={\rm c}_1(\mathcal N_{\ol j}^\vee)\in A^1(Z'),$$
for the first Chern class of the conormal sheaf $\mathcal N_{\ol j}^\vee$
 of $\ol j\colon Z'\to X$.
\item[c)]
The scheme $Z'$ is the coarse space of the stack $\wt Z'$.
The natural decomposition of $(a)$ induces a decomposition $$Z'=\textsqcup_{q=0}^{r-1}
Z_q'.$$
We write  $$\ol j_q\colon Z_q'\to X$$
for the restriction of $\ol j\colon Z'\to X$ to
$Z'_q$.
We write $\tau_q$ for the restriction of $\tau$ to $Z_q$
(when there is no ambiguity we simply write $\tau$).
\item[d)]
Let $\wt L$ be the line bundle over $\wt Z'$ whose fibre
over a point of $p'\in \wt Z'$ is the line
cotangent to the prescribed branch associated to $p'$.
The cotangent lines to the remaining branch yield the line bundle
$\wt{\iota}^* \wt L$. By $\ol L$ and $\ol \iota^* \ol L$,
we denote the analogue line bundles on  $Z'$.
Since $Z'$ is the coarse space of $\wt Z'$, we can pull
$L$ back to $Z'$ via $p\colon \wt Z' \to Z'$; we get
$$p^* \ol L = \wt L^{\otimes r}.$$
We write $$\psi={\rm c}_1(\ol L), \quad \widehat\psi={\rm c}_1(\ol \iota^*\ol L) \quad\text{in } A^1(Z_q'),$$
Clearly, we have
$\widehat\psi = \ol \iota^* \psi$.
 \item[e)]
Consider the above desingularization $\nu\colon C\to \ol C$ (see Figure 1).
The exceptional locus ${\rm Ex}(\nu )\subset C$ is a
curve over $Z\subset \ol C$ whose fibre is
a chain of $r-1$ rational lines.
It can be
regarded as a subcurve of $C\times _X Z\to Z$.

We pullback these $Z$-curves to $Z'$ via $\ol \ep \colon Z'\to Z$.
\[ \xymatrix{
 {\rm Ex}(\nu)\times _Z Z' \ar[r]& C\times _X Z'\ar[d]^{\nu}
\\
& \ol C\times_X Z'\ar[d]\\
& Z'. \ar@/^0.5cm/[u]^{f}\ar@/^1cm/[uul]^{f'}\\
} \]
The curve
$\ol C\times_X Z'\to Z'$ has a distinguished section
$f\colon Z'\to \ol C\times_X Z'$
in the singular locus and, by the definition of $Z'$,
all along this section there is a distinguished
choice of a branch meeting $f(Z')$.
This induces a section $f'$ mapping $Z'$ to the singular locus of $C\times_X Z'$
so that  for any $p'\in Z'$ the node $f'(p')$ lifts $f(p')$ and has a branch lying over
the chosen branch of $f(p')$.

In the above
curve ${\rm Ex}(\nu)\times _Z Z' \longrightarrow Z'$ each fibre is a
chain of $r-1$ rational lines, which can now be numbered
with indices $i$ running from $1$ to $r-1$ starting from the node $f'(Z')$.
We write $$E'_{i}\to Z'$$ for the $\PP^1$-bundles
on $Z'$ obtained from
the normalization of  ${\rm Ex}(\nu)\times _Z Z'\longrightarrow Z'$
(see Figure 2).

\begin{figure}[h]
\begin{picture}(200,200)(50,-30)
\put(163,150){$E_{1}'$}
\put(198,162){\circle*{2}}
\put(168,132){\circle*{2}}
  \qbezier(160,125)(205,170)(205,170)
\put(193,110){$E_{2}'$}
\put(198,132){\circle*{2}}
\put(228,102){\circle*{2}}
\put(228,132){\circle*{2}}
\put(288,102){\circle*{2}}
  \qbezier(190,140)(235,95)(235,95)
\put(223,150){$E_{3}'$}
\put(258,162){\circle*{2}}
  \qbezier(220,125)(265,170)(265,170)
\put(253,110){$E_{4}'$}
\put(258,132){\circle*{2}}
  \qbezier(250,140)(295,95)(295,95)
  \put(290,135){$...$}
\put(332,110){$E_{r-2}'$}
\put(403,162){\circle*{2}}
  \qbezier(365,125)(410,170)(410,170)
\put(302,150){$E_{r-3}'$}
\put(343,162){\circle*{2}}
\put(313,132){\circle*{2}}
\put(343,132){\circle*{2}}
  \qbezier(335,140)(380,95)(380,95)
\put(362,150){$E_{r-1}'$}
\put(373,132){\circle*{2}}
\put(373,102){\circle*{2}}
  \qbezier(305,125)(350,170)(350,170)
\put(280,80){\vector(0,-1){40}}
\put(285,60){\text{normalization}}
  \qbezier[30](395,-15)(405,20)(445,30)
  \qbezier[30](175,-15)(165,20)(125,30)
  \qbezier(160,-15)(205,30)(205,30)
\put(162,-22){\footnotesize$(0)$}
  \qbezier(190,30)(235,-15)(235,-15)
\put(193,8){\footnotesize$(1)$}
  \qbezier(220,-15)(265,30)(265,30)
\put(222,-22){\footnotesize$(2)$}
  \qbezier(250,30)(295,-15)(295,-15)
\put(252,8){\footnotesize$(3)$}
  \put(290,10){$...$}
  \qbezier(410,-15)(365,30)(365,30)
\put(328,-25){\footnotesize$(r-3)$}
  \qbezier(380,30)(335,-15)(335,-15)
\put(340,20){\footnotesize$(r-2)$}
  \qbezier(350,-15)(305,30)(305,30)
\put(397,-22){\footnotesize$(r-1)$}
 \end{picture}
\caption{
The normalization of the fibre over $Z'$:
}
\end{figure}
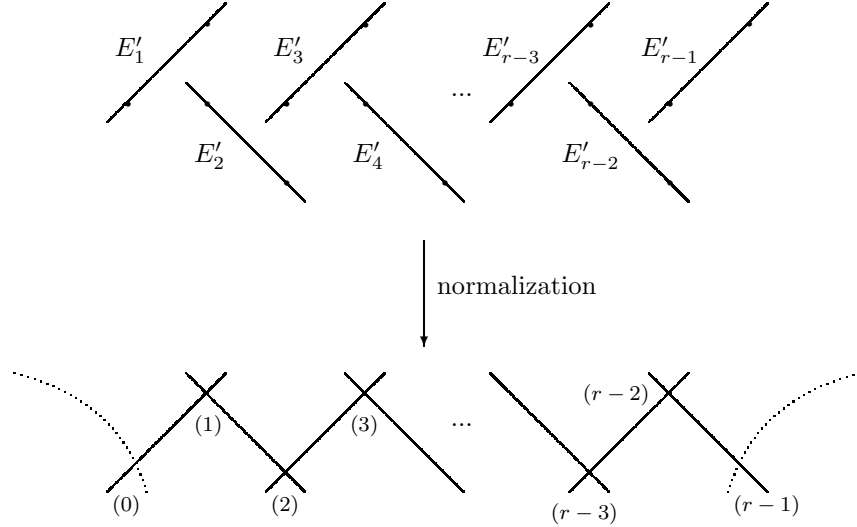

Pulling back $E_i'\to Z'$ via $Z'_q\into Z$ we finally get
$$ E'_{q,i}\to Z'_q.$$

\item[f)]
Note that the fibre of the exceptional locus
${\rm Ex}(\nu)\times _Z Z'$ over $Z'$
contains nodes of $C\times_X Z'$ which are
naturally labelled from $0$ to $r-1$
starting from
$f'(Z')$. In this way the singular locus of $C\times_X Z'\to Z'$ can be regarded
as the disjoint union of $r$ copies $(Z')^{(i)}$ of $Z'$ for
$i\in \{0,\dots, r-1\}$. Over $Z'_q$ we have
$r$ copies $(Z'_q)^{(0)},\dots, (Z'_q)^{(r-1)}$.

We write
\begin{equation}\label{eq:ladef}
\la _{q,i}\colon Z'_q\to E'_{q,i} \quad (\text{ for } i\ne 0),
\quad \quad \text {and} \quad  \quad
\wh \la _{q,i}\colon Z'_q\to E'_{q,i+1} \quad (\text{ for } i\ne r-1)
\end{equation}
for the sections lifting the $i$th node $(Z_q')^{(i)}$
to the bundles $E'_{q,i}$ and $E'_{q,i+1}$
(see Figure 2).
We write $L_{q,i}$ and  $\widehat L_{q,i}$  for the corresponding conormal sheaves.
By setting $L_{q,0}=\ol L$ and $\widehat L_{q,r-1}=\ol \iota ^*\ol L$ we extend the
definition of $L_{q,i}$ and $\widehat L_{q,i}$ to
all $0\le q,i\le r-1$. We set
$$\psi_{q,i}=
{\rm c}_1(L_{q,i}),\quad  \widehat{\psi}_{q,i}=
{\rm c}_1(\widehat L_{q,i})\quad \text{ in } A^1(Z_q'), \quad \text{ for }
0\le q,i\le r-1.$$

\item[g)]
Note that the involution $\iota$ interchanging the
two sheets of $Z'$ yields a natural involution of
${\rm Ex}(\nu)\times _Z Z' \longrightarrow Z'$ and of its normalization. The
involution maps $Z'_q$ to $Z'_{\ol{r-q}}$ and the $i$th component of the
curve ${\rm Ex}(\nu)\times _Z Z' \longrightarrow Z'$
to the $(r-i)$th component.
In this way, the above $\PP^1$-bundle $E'_{q,i}\to Z'_q$
is sent to the $\PP^1$-bundle $E'_{\ol {r- q},r-i}\to Z'_{\ol{r-q}}$.
These two $\PP^1$-bundles are sent via the projection to $C$
to the same divisor, which we denote by
$$E_{q,i}\in \Div (C).$$
Clearly $E_{q,i}$ and $E_{\ol{r-q},r-i}$ denote the same divisor in $C$.
The exceptional locus of the desingularization $C\to \ol  C$
can be written as:
$${\rm Ex}(C\to \ol C)=
\textsum_{q=0}^{\lfloor r/2\rfloor} H_q,$$
where
$$H_q=\begin{cases} \sum_{1\le i\le r-1 } E_{q,i} & \text{for } 0<q<r/2, \\
\\
\sum_{1\le i\le \lfloor r/2\rfloor} E_{q,i} & \text{if $q=0$ or $q=r/2$.}\end{cases}$$

\item[h)] We write $$\sigma_{q,i}\colon Z'_q\to E_{q,i} \quad \quad \text{and}
\quad \quad
\wh \sigma_{q,i}\colon Z'_q\to E_{q,i+1}$$
for the composition
of $\la_{q,i}$  and $\wh \la_{q,i}$ of \eqref{eq:ladef} with the normalization. We write
$$\rho_{q,i}\colon E_{q,i}\into C$$
for the natural inclusion of $E_{q,i}$ in $C$.
Note that $\rho _{q,i}\circ  \sigma_{q,i}$ and
$\rho _{q,i}\circ \wh \sigma_{q,i}$
coincide; we denote such composite morphism  by
$$j_{q,i}\colon Z'_q\to C.$$ The morphisms above fit in the following commutative diagram
$$\xymatrix@C=2cm{
E_{q,i}'\ar[r]\ar[d] & E_{q,i} \ar[r]^{\rho_{q,i}} & C \ar[d]\\
Z_q'\ar[r]\ar[r]\ar@/^0.5cm/[u]^{\la_{q,i}}\ar[ru]^{\sigma_{q,i}}\ar[rru]_{j_{q,i}}
& {Z}_q \ar[r] & \olC}$$
(the diagram remains commutative if we replace
$\lambda_{q,i}, \sigma_{q,i},$ and
$j_{q,i}$ with
$\wh \la_{q,i-1}, \wh \sigma_{q,i-1},$ and
$j_{q,i-1}$).

\item[i)] The characteristic classes $\psi_{q,i}, \wh \psi_{q,i}, \tau$
satisfy the following relations:
\begin{align}
&\widehat \psi_{q,i} = -\psi_{q,i+1},\label{eq:psiopp}\\
&\tau=\psi_{q,i}+\widehat \psi_{q,i}.\label{eq:tauplus}
\end{align}
The first relation follows from the fact that $\psi_{q,i}$ and $\psi_{q,i+1}$
are the first Chern classes of two conormal line bundles corresponding to two disjoint
sections of the $\PP^1$-bundle $E'_{q,i+1}\to Z'_q$. The two line bundles
are dual to each other; so, their respective first Chern classes are opposite.

The second relation also involves
first Chern classes.
On the left hand side we have ${\rm c}_1$ of
the conormal line bundle
$\mathcal N_{\ol j}^\vee$ with respect to $\ol j\colon Z'\to X$.
On the right hand side we have ${\rm c}_1$ of
of $(\det \mathcal N_{f'})^\vee,$ where $f'$ is the morphism $Z'\to C'=C\times_X Z'$.
The two line bundles are isomorphic; indeed,
by the definition of the normal sheaf, we have
\begin{multline}\label{eq:arguxy=t}
\det {\mathcal N_{f'}}\cong (\det T_{Z'})^\vee
\otimes f^* \det T_{C'}\cong
\mathcal N_{\ol j}\otimes ({\ol j}^*\det T _{X})^{\vee} \otimes f^* \det T_{C'}\\\cong
\mathcal N_{\ol j}\otimes f^* ( \det T_{C'} \otimes (\pi^*\det T _{X})^{\vee})\cong \mathcal N_{\ol j}\otimes f^* \omega _{\wt C/X}
\cong \mathcal N_{\ol j},
\end{multline}
where the last equality holds because $f^* \omega _{\wt C/X}$ is trivial
(recall Remark \ref{rem:omega} and the definition of
the relative cotangent sheaf of a nodal curve; the claim extends to
$r$-stable curves by \eqref{eq:omegapull}).

As an immediate  consequence of the above equations,
we have the relation $\psi_{q,i+1}=\psi_{q,i}-\tau,$ for all $i=0,\dots, r-1.$
Therefore, using $\psi_{q,0}=c_1(\ol L)=\psi$, we get
$$\psi_{q,i}=\psi -i \tau.$$
Furthermore, summing up the above equations
\eqref{eq:tauplus} for $i=0,\dots r-1$, we get
$\tau=({\psi +\widehat \psi})/{r}$. Hence we have
\begin{align}
&\psi_{q,i}= \frac{(r-i)\psi - i\widehat{\psi} }{r} \label{eq:psiqi}
\end{align}
\end{enumerate}
\end{rem}

We conclude this section with a result allowing  us to
handle the Grothendieck Riemann--Roch calculation.
We define the \emph{$q$-chain divisor $D_q$} needed in the statement.
We illustrate the definition with a picture. Let us fix $q\in ]1,r/2[$ for simplicity;
then we define $D_q$ by summing the divisors $E_{q,i}$ for $i$ ranging from $r-1$ to $1$
with multiplicities ranging from $q$ to $r-q$ as follows.
We start with $E_{q,r-1}$ with multiplicity $q$ and we increase the multiplicity
with step $q$ up to
$(r-q)q E_{q,q}$; then, we decrease with step $r-q$
until we reach $(r-q) E_{q,1}$.
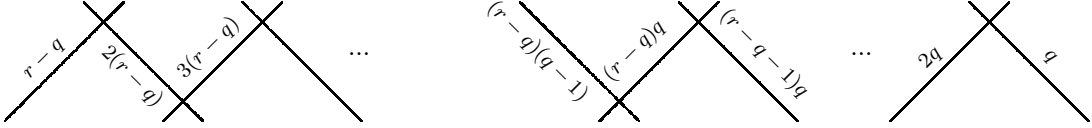
\begin{figure}[h]
\begin{picture}(200,70)(125,0)
  \qbezier(160,5)(205,50)(205,50)
\put(195,31){\rotatebox{315}{\footnotesize$2(r-q)$}}
  \qbezier(190,50)(235,5)(235,5)
\put(166,21){\rotatebox{45}{\footnotesize$r-q$}}
\put(224,21){\rotatebox{45}{\footnotesize$3(r-q)$}}
\put(385,20){\rotatebox{45}{\footnotesize $(r-q)q$}}
\put(341,48){\rotatebox{315}{\footnotesize $(r-q)(q-1)$}}
\put(429,45){\rotatebox{315}{\footnotesize $(r-q-1)q$}}
\put(505,25){\rotatebox{45}{\footnotesize $2q$}}
\put(552,30){\rotatebox{315}{\footnotesize$q$}}
  \qbezier(220,5)(265,50)(265,50)
  \qbezier(250,50)(295,5)(295,5)
  \put(290,30){$...$}
  \qbezier(460,5)(415,50)(415,50)
  \qbezier(430,50)(385,5)(385,5)
  \qbezier(400,5)(355,50)(355,50)
  \put(480,30){$...$}
  \qbezier(540,50)(495,5)(495,5)
  \qbezier(570,5)(525,50)(525,50)
 \end{picture}
\caption{the $q$-chain divisor.}
 \end{figure}

We give  the formal definition
of $D_q$:
\begin{equation}\label{eq:qdiv}
D_q= \begin{cases}
\sum _{1\le i\le q} (r-q)i E_{q,i} + \sum _{q< i <r} q(r-i) E_{q,i} & \text{for }
0<  q<r  \text{ and } q\ne r/2,\\ \\
\sum _{1\le i\le q } (r-q)i E_{q,i}=\sum _{q< i <r} q(r-i) E_{q,i}
  & \text{for $q=r/2$ \text{ and } $q=0$.}\\
\end{cases}\end{equation}
Clearly $D_0$ vanishes by definition, whereas $D_q=D_{r-q}$
for $0<q<r$ and $q\ne r/2$.

\begin{exa}
Let us consider a nodal curve given by the union of two curves
$C'$ and $C''$ meeting at one node.
Assume that the node
is replaced by a chain $H$ of $r-1$ projective lines as above
with the first component $E_{q,1}$ attached to $C'$ and the
last component $E_{q,r-1}$ attached to $C''$.
Then, the $q$-chain divisor supported on $H$
has degree $-r$ on the component $E_{q,q}$,
$0$ on the remaining components of the chain,
$r-q$ on $C'$, and $q$ on $C''$.
\end{exa}
\begin{lem}\label{lem:Dq}
In the situation described above,
let $\wt S$ be the $r$th root on $\wt C$, let $\ol S$
be its direct image on
$\ol C$, and let $\nu\colon C\to \ol C$ be the desingularization.
There exists a line bundle $S$ on $C$ satisfying
\begin{equation}\label{eq:LrootD}
S^{\otimes r}\cong \nu^*\mathcal K \otimes
\Ocal(-\textsum _{1\le q\le \lfloor r/2\rfloor} D_q)
\end{equation}
and
\begin{equation}\label{eq:LolL}
\nu_* S=\ol {S}.
\end{equation}
\end{lem}
\begin{proof}
The statement follows from results of Cornalba \cite{Co} for
square roots generalized to
the case $r>2$ by Jarvis \cite{Ja_geom}, \cite{Ja_Pic} and
by Caporaso, Casagrande, and Cornalba in \cite{CCC}.
We contract the divisors
$$E_{q,i} \text{ for } 0\le q\le  \lfloor r/2\rfloor \text{ and } i\ne q$$
(in particular all divisors $E_{0,i}$ are contracted).
In this way, we obtain a semistable curve $Q\to X$ and a morphism $\nu_Q\colon Q\to \ol C$,
whose
exceptional locus is
a divisor fibred in projective lines; more precisely,
by construction, ${\rm Ex}({\nu_Q})$ is the disjoint union
of the images $I_q$ of
$$E_{q,q} \text{ for } 0< q\le  \lfloor r/2\rfloor$$
via $C\to Q$.
As soon as $r>2$, $I_q$ is not a Cartier divisor:
each fibre of $I_q$ on $\ol C$ contains two singular points, whose
local pictures
in $Q$ are given by $\{xy=t^q\}$ and $\{xy=t^{r-q}\}$ in $\CC^{\dim X+2}$.
On the other hand the divisor
$q(r-q)I_q$ is a Cartier divisor in $Q$ and if we regard $I_q$ as a $\PP^1$-bundle,
on each fibre $F\cong \PP^1\subseteq I_q$ on $\ol C$ we have
$$\deg (q(r-q)I_q)\rest {F}= - \frac{q(r-q)}{q} - \frac{q(r-q)}{r-q}=-r.$$

Jarvis \cite[Thm.~3.3.9]{Ja_Pic} shows that there is a line bundle $S_Q$ on
the semistable  curve $Q$
which is sent to $\ol S$ via direct image to $\ol C$, and
which is equipped with an isomorphism
$$S_Q^{\otimes r}\xrightarrow{\ \sim \ } \nu_Q^*\mathcal K\otimes \Ocal( - \textsum _q {q(r-q)} I_q ),$$
(see \cite{Ja_Pic}, in particular
pt.~1-3 after Lem.~3.3.8, the proof of Thm.~3.3.9, and Fig.~1 therein).
Indeed, $Q$ is realized as ${\rm Proj}({\bigoplus} _{k\ge 0} \Sym^k(\ol S))$
and, with respect to such isomorphism,
$S_Q$ is identified with $\Ocal (1)$, see \cite[\S4.2 and Prop.~4.2.2]{CCC} and
Falting's description of torsion-free sheaves in \cite{Fa} and \cite[\S3.1 and \S3.2]{Ja_comp}.
In fact, it should be also noted that the realization
of the direct image of $S_Q$ as a direct image from the stack $\wt C$ is
due to Abramovich and Jarvis, \cite[Prop.~4.3.1]{AJ}.

It is easy to see that
pulling $q(r-q)I_q$ back to $C$ via $C\to Q$ yields the divisor $D_q$;
therefore, the pullback of $S_Q$ to $C$  yields
the line bundle $S$ satisfying the properties stated above.
\end{proof}

\section{Grothendieck Riemann--Roch for the universal $r$th root}\label{sect:GRR}
We work with the moduli stack
$\MMMbar_{g,n}^{r}$ of $r$th roots of $\mathcal K$ on $r$-stable curves.
Let $\sta C$ be the universal $r$-stable curve, and let $\SSS$ be the universal $r$th root.
\subsection{The proof of the main theorem}
\begin{notn}[tautological classes]\label{notn:taut}
For any $r$-stable curve $\wt C\to X$ with
$n$ markings $x_1, \dots, x_n \colon X\to C$, we write
\begin{equation}\label{eq:defkappapsi}
\kappa^{}_d:=\pi_* ({\rm c}_1(\omega^{\log}_{\wt C/X})^{d+1}) \quad \quad \text{and} \quad
\quad\psi_i={\rm c}_1(x_i^*\omega_{\wt C/X}).\end{equation}
By \cite{AC}, we have
\begin{equation}\label{eq:kappalog}
\kappa_d^{ }=\pi_*({\rm c}_1(\omega_{\wt C/X})^{d+1})+\textsum _{i=1}^n \psi_i^d.
\end{equation}
The cohomology classes and relations
introduced above are usually given in the literature
for the category of stable curves
$\MMMbar_{g,n}$; however they equally hold for $r$-stable curves
by \eqref{eq:omegapull}.
\end{notn}

\begin{notn}[Bernoulli polynomials and numbers]\label{notn:bernoulli}
The Bernoulli polynomials $B_n(x)$ are defined by
\begin{equation}\label{eq:berdef}
\frac{e^{tx}t}{e^t-1}=\sum_{n=0}^{\infty} \frac{B_n(x)}{n!}t^n.
\end{equation}
We recall that the Bernoulli polynomials
satisfy the  property $B_n(x)=(-1)^nB_n(1-x)$ and
\begin{equation}\label{eq:berec}
B_n(x+y)=\sum_{m=0}^n \binom{n}{m}B_m(x)y^{n-m},
\end{equation}
which can be regarded as a formula expressing Bernoulli polynomials
in terms of Bernoulli numbers $B_n(0)$.
\end{notn}

\begin{proof}[Proof of Theorem \ref{thm:GRRcalc}]
We prove the Theorem in the \'etale topology of the stack of $r$th roots of $\Kcal$;
therefore, we work directly with $\wt C$, $\ol C$, and $C$ over a base scheme $X$,
and the $r$th root $\wt S$ on $\wt C$.
In order to calculate $\ch(R^\bullet\wt{\pi}_*\wt{S})$ we proceed in three steps as follows.
\begin{enumerate}
\item We consider the desingularization $\pi \colon C\to X$ of $\ol C \to X$.
By Lemma \ref{lem:Dq}, $R^\bullet\wt{\pi}_*\wt{S}$ equals
$R^\bullet{\pi}_*{S}$, for a line bundle $S\in \Pic (C)$ whose
first Chern class can be
expressed in the rational Chow ring in terms of $K={\rm c}_1(\omega_{C/X})$,
the divisors $\Delta_1, \dots, \Delta_n$ specifying the markings,
and the $q$-chain divisors $D_q$.
The calculation boils down to two distinct terms
on which we focus in the second and third step of the proof.
\item We focus on the terms involving $K$ and $\Delta_i$.
\item We conclude the calculation on the boundary locus.
\end{enumerate}

\vspace{0.2cm}

\noindent Step 1.
The curve $C\to X$ is a flat and proper morphism
between regular schemes. Therefore, we can
apply Grothendieck Riemann--Roch (GRR) in the form
$$\ch(R^\bullet \pi_*S)=\pi_*(\ch(S)\td^\vee(\Omega^1_{C/X})).$$
We adopt Mumford's notation for the Todd character:
for any line bundle $R$
we have
$$\td^\vee (R)= \frac{{\rm c}_1(R)}{\exp({{\rm c}_1(R)})-1}.$$
Recall the exact sequence
$0\to\Omega^1_{C/X}\to \omega_{C/X} \to \omega_{C/X}\otimes
\Ocal_{Z}\to 0$ and write
$$\ch(R^\bullet \pi_* S)=\pi_*\left ((\ch S)(\td^\vee\omega_{C/X})(\td^\vee\Ocal_{\Sing(C/X)})^{-1}\right).$$
By \eqref{eq:LrootD}, the right hand side equals
\begin{equation}\label{eq:GRRexplicit}
\pi_*\left(\exp \left(sK/r-\sum_{i=1}^n (m_i-s)\Delta_i /r -\sum_{q=1}^{\lfloor r/2\rfloor} D_q/r\right)
\frac{K}{\exp(K)-1}(\td^\vee\Ocal_{\Sing(C/X)})^{-1}\right),
\end{equation}
where $K$ is ${\rm c}_1(\omega_{C/X})$.
We make the following observations.
\begin{enumerate}
\item For each $i$ the divisor $\Delta_i$ is disjoint from  ${\Sing(C/X)}$.
\item The intersections between $K$ and any cycle in $\Sing(C/X)$ vanish (see Remark \ref{rem:omega}).
\item The intersections between $K$ and any divisor $D_q$ vanish
(in other words $C\to \olC$ is crepant).
\item The divisor $D_q$ is disjoint from $D_{q'}$
for distinct indices $q$ and $q'$ satisfying $0\le q, q'\le \lfloor r/2\rfloor$.
\end{enumerate}
In this way, the term \eqref{eq:GRRexplicit} is the direct image via $\pi_*$ of
\begin{equation}\label{eq:twosummands}
\left(\frac{K\exp (sK/r)}{\exp(K)-1}\exp(-\textsum_{i=1}^n (m_i-s)\Delta_i /r)\right)
+
\left(\left(\exp(-\textsum_{q=1}^{\lfloor r/2\rfloor} D_q/r)(\td^\vee\Ocal_{\Sing(C/X)})^{-1}\right)-1\right).
\end{equation}
Step 2 and Step 3 are devoted to the first and the second summand of \eqref{eq:twosummands},
respectively.

\begin{center}-----------------------------------
\end{center}

\noindent Step 2. We consider the term
\begin{equation}\label{eq:termstep2}
\pi_*\left(\frac{K\exp (sK/r)}{\exp(K)-1}\exp(-\textsum_{i=1}^n (m_i-s)\Delta_i /r)\right).
\end{equation}
Using Bernoulli polynomials we expand the term inside the brackets in
the formula above as the product of
$$\sum_{d\ge -1} \frac{B_{d+1}(\frac{s}{r})}{(d+1)!}K^{d+1} \quad \quad
\text{ and } \quad
\quad
\prod_{i=1}^n \sum _{l\ge 0} \frac{ (m_i-s)^l (-\Delta_i)^l}{l! r^l}.$$
The direct image via $\pi$ of terms of degree $0$ vanishes and the
divisors $\Delta_i$ are disjoint; therefore,
we are only interested in the direct image of the following sum,
$$\sum_{d\ge 0}\frac{B_{d+1}\left(\frac{s}{r}\right)}{(d+1)!}K^{d+1} +
\sum_{i=1}^n\sum _{h,l\ge 0}
\frac{ \big(\frac{m_i-s}{r}\big)^{l+1} B_h\left (\frac{s}{r}\right)(-\Delta_i)^{l+1}K^h}{(l+1)! h!}
$$

By Remark \ref{rem:omega}, for all $i=1, \dots, n$
we have $-K\Delta_i=\Delta_i^2$. Therefore,
$(-\Delta_i)^{l+1}K^h$ equals $-K^{h+l}\Delta_i$
and we rewrite the above sum as
\begin{equation}\label{eq:twosumstep2}\sum_{d\ge 0}
\frac{B_{d+1}\big(\frac{s}{r}\big)}{(d+1)!} K^{d+1} -
\sum_{i=1}^n\sum_{h, l\ge 0}
\frac{ \left(\frac{m_i-s}{r}\right)^{l+1}
B_h\left(\frac{s}{r}\right)\Delta_iK^{h+l}}{(l+1)! h!}.
\end{equation}
We now identify the term of degree $d+1$ in the
above sum (since $\pi $ is of relative dimension $1$
the terms of degree $d+1$ contribute to the degree-$d$ term of $\ch$).
The first summand contains only one term of degree $d+1$: the
class $B_{d+1}(s/r)K^{d+1}/(d+1)!$.
The contribution of the second summand in degree $d$ occurs for each $i=1,\dots, n$ and
equals $K^{d} \Delta_i$ times
\begin{align*}
-&\sum _{h=0}^d
\frac{\big(\frac{m_i-s}{r}\big)^{d-h+1} B_h(\frac{s}{r})}{(d-h+1)! h!}
&=-\frac{1}{(d+1)!}\sum _{h=0}^{d+1}
\textstyle{\binom{d+1}{h} B_h({s}/{r}) \left(\frac{m_i-s}{r}\right)^{d+1-h}+
{B_{d+1}({s}/{r})}
}\end{align*}
Using \eqref{eq:berec} we get
$\textstyle{\sum _{h=0}^{d+1}
\binom{d+1}{h} B_h({s}/{r})\left(\frac{m_i-s}{r}\right)^{d+1-h}=B_h\left(\frac{m_i-s}{r}+\frac{s}{r}\right)=
B_h\left(\frac{m_i}{r}\right)}$
and we simplify the above term:
$$\left(-\frac{B_{d+1}\left(\frac{m_i}{r}\right)}{(d+1)!}+\frac{B_{d+1}\left(\frac{s}{r}\right)}{(d+1)!}\right) K^d\Delta_i.$$
By \eqref{eq:kappalog}, $\kappa^{}_d=
\pi_* (K^{d+1})+\sum_{i=1}^n \psi_i^d$.
Recall that $\psi_i^d=\pi_*(K^d\Delta_i)$.
Putting everything together we express the term \eqref{eq:termstep2} as
$$\frac{B_{d+1}\left(\frac{s}{r}\right)}{(d+1)!}
\kappa_d^{} -\sum_{i=1}^n\frac{B_{d+1}\big(\frac{m_i}{r}\big)}{(d+1)!} \psi_i^{d}.$$

\begin{center}-----------------------------------
\end{center}

\noindent Step 3.
We now consider the second summand in \eqref{eq:twosummands}. We omit $\pi_*(1)$
because $\pi$ is a morphism of relative dimension $1$. Therefore, we
focus on the degree-$d$ term of
\begin{equation}\label{eq:lastterm}
\pi_*
\left(\exp(-\textsum_{1\le q\le \lfloor r/2\rfloor} D_q/r)(\td^\vee\Ocal_{\Sing(C/X)})^{-1}\right)
\end{equation}
and show that it equals
\begin{equation}\label{eq:lasttermtarget}
\frac{1}{2}\sum_{q=0}^{r-1}\frac{rB_{d+1}(q/r)}{(d+1)!}(\ol j_q)_*(\gamma_{d-1})
\end{equation}
By Remark \ref{notn:charclsing}.(f)
we regard the singular locus of $C\times _X Z'$ as the union of
copies
$(Z_q')^{(i)}$ of $Z_q'$ for $i=0,\dots , r-1$. Then, we
can rewrite Mumford's formula \cite[Lem.~5.1]{Mu_GRR} for
$ (\td^\vee\Ocal_{\Sing(C/X)})^{-1}$
as
\begin{equation}\label{eq:MuSing}
(\td^\vee\Ocal_{\Sing(C/X)})^{-1}=1+\frac{1}{2}\sum_{m\ge 1} \ \sum_{0\le q,i\le r-1} \frac{B_{m+1}(0)}{(m+1)!}
(j_{q,i})_*\left(\frac{\psi_{q,i}^m+\widehat
{\psi}_{q,i}^m}{\psi_{q,i}+\widehat\psi_{q,i}}\right),\end{equation}
(see Remark \ref{notn:charclsing}.(h) for the definition of $j_{q,i}$).
Then,
we replace $(\td^\vee\Ocal_{\Sing(C/X)})^{-1}$ in \eqref{eq:lastterm}
$$\pi_*
\left(\left(\prod_{q=1}^{\lfloor r/2\rfloor}\sum_{h\ge 0} \frac{(-D_q/r)^h}{h!}
\right)\left(1+\frac{1}{2}\sum _{0\le q,i<r}
\ \sum_{m\ge 1} \frac{B_{m+1}(0)}{(m+1)!}(j_{q,i})_*\left(\frac{\psi_{q,i}^m+\widehat
{\psi}_{q,i}^m}{\psi_{q,i}+\widehat\psi_{q,i}}\right)\right)\right)$$
and we expand using the fact that the divisors $D_1, \dots, D_{\lfloor r/2\rfloor}$
are disjoint
\begin{equation}\label{eq:beforelem1}\pi_*
\left(\left(1 + \sum_{q=1}^{\lfloor r/2\rfloor} \sum_{h\ge 0} \frac{(-D_q/r)^{h+1}}{(h+1)!}
\right)
\left(1+\frac{1}{2}\sum _{0\le q,i<r}
\ \sum_{m\ge 1} \frac{B_{m+1}(0)}{(m+1)!}(j_{q,i})_*\left(\frac{\psi_{q,i}^m+\widehat
{\psi}_{q,i}^m}{\psi_{q,i}+\widehat\psi_{q,i}}\right)\right)\right).
\end{equation}
We adopt the convention \begin{equation}\label{eq:convqhat}
\qquad \wh q=\ol{r-q} \qquad \text{for all $q\in \{0, \dots, r-1\}$}.
\end{equation} In this way, we have the following lemma.
\begin{lem}\label{lem:qqcheck} For $0\le q,i\le r$
we have
$$(j_{q,i})^*(-D_q)=\begin{cases} \wh q \psi  & \text{for $i< q$},\\
                               q\widehat \psi  &\text {for $i\ge q$}. \end{cases}$$
\end{lem}
\begin{proof}
We use the notation established in Remark \ref{notn:charclsing}.(h) and we calculate $(j_{q,i})^*E_{q,i}$  by factoring
$j_{q,i}$ as the composite of $\rho_{q,i}$ and $\sigma_{q,i}$.
Similarly, we calculate
$(j_{q,i})^*E_{q,i+1}$ by factoring
$j_{q,i}$ as the composite of $\rho_{q,i+1}$ and $\wh \sigma_{q,i}$.
We get
$$(j_{q,i})^*E_{q,i}=(\sigma_{q,i})^*\rho_{q,i}^*E_{q,i}=
(\sigma_{q,i})^*\mathcal {\rm c}_1(\mathcal N_{E_{q,i}/C})$$
and, similarly,  $(j_{q,i})^*E_{q,i+1}=(\wh{\sigma}_{q,i})^*
{\rm c}_1(\mathcal N_{E_{q,i+1}/C}).$
Therefore,
$(j_{q,i})^*E_{q,i}$ is the first Chern class of the line bundle whose fibres are
tangent along the branches orthogonal to $E_{q,i}$. In this way, we have
$$(j_{q,i})^*E_{q,i}=-\widehat \psi_{q,i} \quad \text{and} \quad
(j_{q,i})^*E_{q,i+1}=-\psi_{q,i}.$$
This implies the desired realtions.
Consider the case $1\le q< r/2$.
We apply the definition of $D_q$ and the
relations of Remark \ref{notn:charclsing}.(f); we get
\begin{align*}
-(j_{q,i})^*D_q&=\textsum_{i=1}^{q} (r-q)i (j_{q,i})^*(-E_{q,i}) +\textsum_{i=q+1}^{r-1} q(r-i) (j_{q,i})^*(-E_{q,i})\\
&=\begin{cases}(r-q)(i \widehat{\psi}_{q,i}+(i+1)\psi_{q,i}) & \text{if $i\le q-1$}\\
q((r-i) \widehat{\psi}_{q,i}+(r-i-1)\psi_{q,i}) & \text{if $i\ge q$}\end{cases}\\
&=\begin{cases}(r-q)(i \tau+ \psi_{q,i}) & \text{if $i\le q-1$}\\
q((r-i) \tau- \psi_{q,i}) & \text{if $i\ge q$}\end{cases}\\
&=\begin{cases}(r-q)\psi\rest{\pri Z_q} & \text{if $i\le q-1$}\\
q (r\tau- \psi)\rest{\pri Z_q}  & \text{if $i\ge q$}\end{cases}.
\end{align*}
\end{proof}
We expand the formula obtained in \eqref{eq:beforelem1}
\begin{multline}
\pi_*\Bigg(
1 +
\sum_{\substack{1\le q\le \lfloor r/2\rfloor\\ d\ge0 }} \frac{(-D_q/r)^{d+1}}{(d+1)!}
\\+\frac{1}{2}\sum_{\substack{0\le q, i<r}}\
\sum_{\substack{m\ge 1\\ h\ge 0}} \frac{B_{m+1}(0)}{(m+1)!h!}
(j_{q,i})_*\Bigg(\sum_{0\le i <r}\frac{\psi_{q,i}^m+\widehat
{\psi}_{q,i}^m}{\psi_{q,i}+\widehat\psi_{q,i}}(j_{q,i})^*(-D_{q}/r)^h\Bigg)
\Bigg).
\end{multline}
We split the summation over $i$ into a first summation ranging
over the indices
$i<q$ and a second summation ranging over the indices $i\ge q$.
In this way, we can apply Lemma \ref{lem:qqcheck} and
the above equation can be rewritten as
follows (we omit $\pi_*(1)$ because it vanishes and we
write $\pi\circ j_{q,i}$ as $\ol j_q$)
\begin{multline}
\pi_*\Bigg(
\sum_{\substack{1\le q\le \lfloor r/2\rfloor\\ d\ge0 }}
\frac{(-D_q/r)^{d+1}}{(d+1)!}\Bigg)
\\+ \frac{1}{2(d+1)!}
\sum_{\substack{0\le q,i<r}}(\ol j_{q})_*\Bigg(\sum_{m= 1}^d
{\textstyle{B_{m+1}(0)}\binom{d+1}{m+1}}\sum_{i=0}^{q-1}\frac{\psi_{q,i}^m+\widehat
{\psi}_{q,i}^m}{\psi_{q,i}+\widehat\psi_{q,i}}(\wh q \psi/r )^{d-m}+\\
\sum_{m= 1}^d{\textstyle{B_{m+1}(0)}\binom{d+1}{m+1}}\sum_{i=q}^{r-1}\frac{\psi_{q,i}^m+\widehat
{\psi}_{q,i}^m}{\psi_{q,i}+\widehat\psi_{q,i}}(q \wh \psi/r)^{d-m}\Bigg).
\end{multline}
Note that, when $m$ is odd, $\wh \psi_{q,i}^m$ is opposite to
$\psi_{q,i+1}^m$ and recall that $\psi_{q,i}+\wh \psi_{q,i}$ equals
$\tau=\psi/r+\wh \psi/r$.
Therefore the sum of the terms $\psi_{q,i}^m+\widehat
{\psi}_{q,i}^m$ over the indices $m\in \{0,\dots, q-1\}$ yields
$\psi_{q,0}^m+\widehat
{\psi}_{q,q-1}^m$, hence $\psi_{q,0}^m-
{\psi}_{q,q}^m$.
In this way, using the relation
$\psi_{q,q }=\wh q \psi/r -q\wh \psi/r$ from
\eqref{eq:psiqi}, we rewrite the term between brackets in the
second summand above as the sum $\alpha_{d-1,q}+\wh \alpha_{d-1,q}$,
where $\al_{d,q}$ is defined as
$$\alpha_{d,q}=\frac{\sum_{m= 0}^{d}{\textstyle{B_{m+2}(0)}\binom{d+2}{m+2}}
((r\psi)^{m+1} -(\wh q \psi -q\wh \psi)^{m+1}) (\wh q \psi)^{d-m}}{r^{d}(\psi+\wh\psi)}$$
and
$\wh \alpha_{d,q}$ is defined as $\al_{d,q}$ after the exchanges
$\psi\leftrightarrow \wh \psi$ and $q\leftrightarrow \wh q$.
Finally, we simplify the expression for $\al_{d,q}$ using $a^{m+1}-b^{m+1}=(a-b)\sum_{i=1}^m a^i b^{m-i}$
and $(r\psi) -(\wh q \psi -q\wh \psi)=q\psi+q\wh \psi$ :
\begin{equation}\label{eq:alpha}
\al_{d,q}=
\sum_{m=0}^{d}\frac{{\textstyle{B_{m+2}(0)}}}{r^{d}}\binom{d+2}{m+2}q(\wh q \psi)^{d-m}
\sum_{i+j=m} (r\psi)^i(\wh q \psi -q\wh \psi)^{j}.
\end{equation}
\begin{lem}
For $d\ge 1$,  we have
$$\sum _{q=1}^{\lfloor r/2\rfloor} {\pi_*((-D_q/r)^{d+1})}=
-\frac{1}{2}\sum_{q=0}^{r-1}(\ol j_q)_* \left(\frac{q \wh q}{r^d}
\sum_{a+b=d-1} (\wh q \psi)^{a}(q \wh \psi)^{b}\right).$$
\end{lem}
\begin{proof}
The proof uses the morphisms $\rho_{q,i}$, $\sigma_{q,i}$, and $\wh \sigma_{q,i-1}$
introduced in
Remark \ref{notn:charclsing}.(h). We recall the commutative diagram
\[\xymatrix@C=1.5cm{
&E_{q,i}\ar[d]^\eta \ar[r]^{\rho_{q,i}} &C\ar[d]_\nu \ar@/^0.5cm/[dd]^{\pi} \\
Z_q'\ar@/_.5cm/[drr]_{\ol j_q}
\ar@/^0.8cm/[ur]^{\sigma_{q,i}}\ar@/^0.4cm/[ur]_{\wh \sigma_{q,i-1}}
\ar[r]_{\ol \ep= \varepsilon_{\mid Z_q'}} &Z\ar[dr]_{\ol j_Z} \ar[r]&\olC\ar[d]_{\overline{\pi}}\\
&&X,
}\]
where $\eta$ denotes  the projection of
the $\PP^1$-bundle $E_{q,i}$ to $Z$ and $\ol j_Z$ is the restriction of
$\ol \pi$ to $Z\subset \ol C$.

We write the cycle $E_{q,i}$ in $C$ as
$$E_{q,i}=\pi^*(\ol j_Z)_*(1)- \textsum _{j\ne i} E_{q,j}$$
and we calculate the restriction of $E_{q,i}$ to itself:
\begin{multline*}
({\rho_{q,i}})^*E_{q,i} = ({\rho_{q,i}})^*\pi^*(\ol j_Z)_*(1)- ({\rho_{q,i}})^* \textsum _{j\ne i} E_{q,j}
= \eta^*(\ol j_Z)^*(\ol j_Z)_*(1)- ({\rho_{q,i}})^* \textsum _{j\ne i} E_{q,j}\\
=- \eta^*\tau_Z - \begin{cases}(\wh \sigma_{q,i-1})_*(1)& q=i=r/2,\\
 (\wh \sigma_{q,i-1})_*(1)+ (\sigma_{q,i})_*(1)& \text{otherwise,}\end{cases}
\end{multline*}
where $\tau_Z=-(\ol j_Z)^*(\ol j_Z)_*(1)$ is the first Chern class of the conormal
bundle relative to  $\ol j_Z\colon Z\to X$.
We now consider the pullback of  $D_q$; let $1\le q\le r/2$; the definition of $D_q$ yields
\begin{equation*}
(\rho_{q,i})^*D_q=\begin{cases}
(r-q)(-i \eta^*\tau_Z  - (\wh \sigma_{q,i-1})_*(1)+ (\sigma_{q,i})_*(1))&\text{$1\le i< q$,}\\
-(r-q)q \eta^*\tau_Z  - (r-q)(\wh \sigma_{q,q-1})_*(1)- q(\sigma_{q,q})_*(1)& \text{$i=q\ne r/2$,}\\
-(r-q)q \eta^*\tau_Z  - (r-q)(\wh\sigma_{q,q-1})_*(1)& \text{$i=q=r/2$,}\\
q(-(r-i)\eta^*\tau_Z  + (\wh \sigma_{q,i-1})_*(1)- (\sigma_{q,i})_*(1))& \text{$q<i< r$.}
\end{cases}
\end{equation*}
For $1\le q<r/2$, we have
\begin{align}\label{eq:forDq} (-D_q)^{d+1}=& -\sum_{i=1}^q (r-q)i (\rho_{q,i})_*((\rho_{q,i})^*(-D_q)^d)
-\sum_{i=q+1}^{r-1} q(r-i) (\rho_{q,i})_*((\rho_{q,i})^*(-D_q)^d),
\end{align}
which can be written explicitly as
\begin{multline*}
-\sum_{i=1}^{q-1} (r-q)^{d+1} i (\rho_{q,i})_*\big((\wh\sigma_{q,i-1})_*(1)- (\sigma_{q,i})_*(1)+i
\eta^*\tau_Z \big)^d\\
- (r-q)q(\rho_{q,q})_*\big((r-q)(\wh\sigma_{q,q-1})_*(1)+q(\sigma_{q,q})_*(1)+(r-q)q \eta^*\tau_Z\big)^d\\
-\sum_{i=q+1}^{r-1} q^{d+1} (r-i) (\rho_{q,i})_*\big((\sigma_{q,i})_*(1)-(\wh \sigma_{q,i-1})_*(1) +
(r-i)\eta^*\tau_Z \big)^d
\end{multline*}
(for $q=r/2$, the above formula holds if we erase the third line and the term $q(\sigma_{q,q})_*(1)$
in the second line).
We finally expand each summand in the equation above. Recall that, if $h$ is positive
and $a, b,$ and $c$ commute and satisfy $bc=0$, we have
\begin{equation}\label{eq:comm}
(a+b+c)^d=a^d+\sum_{m=0}^{d-1}\binom{d}{m} a^m (b^{d-m}+c^{d-m}).
\end{equation}
For $t>0$, we have the relations
$$((\wh \sigma_{q,i-1})_*(1))^t=(\wh\sigma_{q,i-1})_*((-\widehat{\psi}_{q,i-1})^{t-1}) \qquad
\text{and}\qquad
(-(\sigma_{q,i})_*(1))^t=-(\sigma_{q,i})_*((\psi_{q,i})^{t-1}).
$$
So, for any $t>0$ and $l\ge 0$, we have
$$
(\eta^*\tau_Z)^l(-(\sigma_{q,i})_*(1))^t=-(\sigma_{q,i})_*(({\psi}_{q,i})^{t-1}(\sigma_{q,i})^*(\eta^*\tau_Z)^l)=
-(\sigma_{q,i})_*(({\psi}_{q,i})^{t-1}\tau^l).$$
and similarly for $\wh \sigma_{q,i-1}$.
We are interested in the direct image of $(-D_q)^{d+1}$ via $\pi$; we remark that
$\pi_*\eta^*\tau_Z=(\ol j_Z)_*\eta_*\eta^*\tau_Z$ vanishes
(as a consequence of $\eta_*\eta^*=0$).
We also note that both $\pi_*(\rho_{q,q})_* (\wh\sigma_{q,i-1})_*$ and
$\pi_*(\rho_{q,q})_* (\wh \sigma_{q,i})_*$
equal $(\ol j_q)_*$.
Finally, we recall the relation
$\psi_{q,i}=-\widehat{\psi}_{q,i-1}$ from \eqref{eq:psiopp},
and we deduce that
the direct image via  $\pi$ of all summands in the
above formulae \eqref{eq:forDq}
corresponding to $i\ne q$ vanish
(the key fact is $\wh \psi_{q,i-1}+\psi_{q,i}=0$, which implies
$$(\ol j_q)_*\left(
(i\tau)^m \left((-\widehat{\psi}_{q,i-1})^{d-1-m} -  \psi_{q,i} ^{d-1-m} \right)\right)=0 \quad \quad
\text{for $0\le m\le d-1$}). $$

The direct image $\pi_*((-D_q)^{d+1})$ boils down to
$$\begin{cases}
- (r-q)q\pi_*(\rho_{q,q})_*\big((r-q)(\wh\sigma_{q,q-1})_*(1)+q(\sigma_{q,q})_*(1)+(r-q)q \tau\big)^d
& q< r/2 \\
- (r-q)q\pi_*(\rho_{q,q})_*\big((r-q)(\wh\sigma_{q,q-1})_*(1)+(r-q)q \tau\big)^d &
q=r/2
\end{cases}
$$
Then \eqref{eq:comm} yields
\begin{multline}
\pi_*((-D_q)^{d+1})= \begin{cases}
- (r-q)^{d+1}q \sum _{m=0}^{d-1} \binom {d}{m} (\ol j_q)_*\left(
(-\widehat{\psi}_{q,q-1})^{d-1-m}(q\tau)^m \right) \\
- q^{d+1}(r-q) \sum _{m=0}^{d-1} \binom {d}{m} (\ol j_q)_*\left(
(-\psi_{q,q}) ^{h-1-m}    ((r-q)\tau )^m\right)  & q<r/2 \\
\\
- (r-q)^{d+1}q \sum _{m=0}^{d-1} \binom {d}{m}(\ol j_q)_*\left(
(-\widehat{\psi}_{q,q-1})^{d-1-m}(q\tau)^m\right) & q=r/2.
\end{cases}
\end{multline}
Therefore, for $\wh q=\ol{r-q}$, we have
$$\sum _{q=1}^{\lfloor r/2\rfloor} {\pi_*((-D_q/r)^{d+1})}=-\frac{1}{2}
\sum_{q=0}^{r-1}(\ol j_q)_* (\beta_{d-1,q}+\wh \beta_{d-1,q})$$
for
\begin{multline*}
\beta_{d,q}= \frac{q \wh q}{r^{d+2}} \sum_{m=0}^{d}\binom {d+1}{m}
\wh q^{d+1} {\psi}_{q,q}^{d-m}(q\tau)^m\\
=\frac{q \wh q}{r^{d+2}} \bigg(\sum_{m=0}^{d+1}\binom {d+1}{m}
\wh q^{d+1} {\psi}_{q,q}^{d-m}(q\tau)^m -\wh q^{d+1} {\psi}_{q,q}^{-1}(q\tau)^{d+1}\bigg)
=\frac{q \wh q}{r^{d+1}} \frac{(\wh q \psi)^{d+1}-((\psi +\wh\psi)q \wh q/r)^{d+1}}{\wh q\psi -q\wh \psi}
\end{multline*}
and for $\wh \beta_{d,q}$ defined as $\beta_{d,q}$ after the usual
exchanges $q\leftrightarrow \wh q$ and $\psi\leftrightarrow \wh\psi$.
We observe that $-((\psi +\wh\psi)q \wh q/r)^{d+1}/({\wh q\psi -q\wh \psi})$ is
antisymmetric with respect to such transformation. So we get
\begin{equation}
\label{eq:beta}\beta_{d,q}+\wh \beta_{d,q}=
\frac{q \wh q}{r^{d+1}} \frac{(\wh q \psi)^{d+1}+(q \wh \psi)^{d+1}}{\wh q\psi -q\wh \psi}=
\frac{q \wh q}{r^{d+1}} \sum_{a+b=d} (\wh q \psi)^{a}(q \wh \psi)^{b}\end{equation}
and the claim.
\end{proof}

We finish the proof by identifying  the degree-$d$ term of \eqref{eq:lastterm} with
the expression \eqref{eq:lasttermtarget}.
We have shown that the term of degree $d$  of \eqref{eq:lastterm} equals
$$\frac{1}{2(d+1)!} \sum_{q=0}^{r-1} \left(\alpha_{d,q}-\beta_{d,q}-\wh \beta_{d,q}+\wh \alpha_{d,q}\right),$$
with $\alpha_{d,q}$ and $\beta_{d,q}$
defined as
in \eqref{eq:alpha} and \eqref{eq:beta}.
By writing explicitly, we get
\begin{multline}\label{eq:thething}\sum_{m=0}^{d}\frac{B_{m+2}(0)}{r^d}\binom{d+2}{m+2}
q(\wh q \psi)^{d-m} \sum_{i+j=m}
(r\psi)^i(\wh q \psi-q\wh \psi)^j \\-\frac{q\wh q}{r^{d+1}}\sum_{i+j=d }
(\wh q\psi )^i(q \wh \psi)^j \\+
\sum_{m=0}^{d}\frac{B_{m+2}(0)}{r^d}\binom{d+2}{m+2}
\wh q(q \wh \psi)^{d-m} \sum_{i+j=m}
(r\wh \psi)^i( q \wh \psi-\wh q \psi)^j, \end{multline}
which can be easily rewritten as $1/(\psi +\wh \psi)$ times
\begin{multline}\label{eq:thething2}
r\sum_{m=0}^{d}{B_{m+2}(0)}\binom{d+2}{m+2}
 \left(\Big(\frac{\wh q \psi}{r}\Big)^{d+1}\Big(\frac{r}{\wh q}\Big)^{m+1}
 - \Big(\frac{\wh q \psi}{r}-\frac{q\wh \psi}{r}\Big)^{m+1}\Big(\frac{\wh q \psi}{r}\Big)^{d-m}\right)
\\-{q}\sum_{i+j=d }
\Big (\frac{\wh q \psi}{r} \Big)^{i+1} \Big(\frac{q\wh \psi}{r}\Big)^j
-{\wh q}\sum_{i+j=d }
\Big(\frac{\wh q \psi}{r}\Big)^{i} \Big(\frac{q\wh \psi}{r} \Big)^{j+1} \\+
r\sum_{m=0}^{d}{B_{m+2}(0)}\binom{d+2}{m+2}
\left( \Big(\frac{q\wh \psi}{r} \Big)^{d+1} \Big(\frac{r}{q}\Big)^{m+1}
- \Big( \frac{q\wh \psi}{r}-\frac{\wh q \psi}{r}\Big)^{m+1}\Big(\frac{q \wh\psi}{r}\Big)^{d-m}\right).
\end{multline}
We regard the above sum as a homogeneous
polynomial of degree $d+1$ in the variables $\wh q\psi/r$ and $q\wh \psi/r$.
We show that, for any $l=1, \dots, d$, the coefficient of the term
$(\wh q\psi/r)^l (q \wh \psi/r)^{d+1-l}$ vanishes. Indeed, it equals
$$r(-1)^{d-l}\sum_{m=d-l}^d \binom{m+1}{d+1-l} \binom{d+2}{m+2}B_{m+2}(0)
-r- r(-1)^l\sum _{m=l-1}^d  \binom{m+1}{l}\binom{d+2}{m+2}B_{m+2}(0),$$
where the three summands correspond to the three lines in the above formula.
We set $a=l$, $b=d+1-l$, and $h=m+1$, and we rewrite the above expression
as
$$-r(-1)^{b}\sum_{h=b}^{a+b} \binom{h}{b} \binom{a+b+1}{h+1}B_{h+1}(0)
-r- r(-1)^a\sum _{h=a}^{a+b}  \binom{h}{a}\binom{a+b+1}{h+1}B_{h+1}(0).$$
Note that $\{(d,l)\mid 0<l\le d\}$
corresponds to $\{(a,b)\mid a,b>0\}$; therefore, we need to
show that the above expression vanishes for all positive integers $a$ and $b$.
The cases $a=1$, $b>0$ follow from \eqref{eq:berec}.
We can conclude by induction,
because the difference between
the above expression for $(a,b)$ and for $(a+1,b-1)$
equals $r$ times
$$(-1)^b\sum_{h=b}^{a+b+1} \binom{h}{b}\binom{a+b+1}{h} B_{h}(0)
+(-1)^a\sum _{h=a+1}^{a+b+1} \binom{h}{a+1}\binom{a+b+1}{h}B_{h}(0),$$
which in turn equals $\binom{a+b+1}{b}$ times
$$(-1)^b\sum_{h=a+1}^{a+b+1} \binom{a+1}{h-b} B_{h}(0)
+(-1)^a\sum _{h=b}^{a+b+1} \binom{b}{h-a-1}B_{h}(0).$$
The above expression  vanishes: for all nonnegative $\al$ and $\be$ the sum
$\label{eq:newber}
(-1)^\al \sum _{i=0}^{\al}\binom{\al}{i} B_{i+\be}(0)$
is
equal to the sum $(-1)^\be \sum _{i=0}^{\be}\binom{\be}{i} B_{i+\al}(0)$
(Carlitz identity).

In order to conclude the proof
it remains to show that
the coefficient of $\psi^{d+1}$ and $\wh \psi^{d+1}$ equal $rB_{d+2}(\wh q/r)$
and $rB_{d+2}( q/r)$, respectively.
Indeed, the coefficient of $\psi^{d+1}$ is equal to
\begin{multline*}
r\sum _{m=0}^d  B_{m+2}(0){\textstyle\binom{d+2}{m+2}}\Big(\frac{\wh q}{r}\Big)^{d-m}
-r\sum _{m=0}^d B_{m+2}(0){\textstyle\binom{d+2}{m+2}}\Big(\frac{\wh q}{r}\Big)^{d+1}
-r B_{d+2}(0)\Big(-\frac{\wh q}{r}\Big)^{d+1}
-q\Big(\frac{\wh q}{r}\Big)^{d+1}=\\
r\sum _{m=0}^d  B_{m+2}(0){\textstyle\binom{d+2}{m+2}}\Big(\frac{\wh q}{r}\Big)^{d-m}
-r\sum _{m=0}^d B_{m+2}(0){\textstyle\binom{d+2}{m+2}}\Big(\frac{\wh q}{r}\Big)^{d+1}
+r B_{d+2}(1)\Big(\frac{\wh q}{r}\Big)^{d+1}
-r\Big(\frac{\wh q}{r}\Big)^{d+1}
+r \Big(\frac{\wh q}{r}\Big)^{d+2}.
\end{multline*}
Applying \eqref{eq:berec} to the three middle terms, we get
$$r\sum _{m=0}^d  B_{m+2}(0)\binom{d+2}{m+2}\left(\frac{\wh q}{r}\right)^{d-m}+
r B_{1}(0)(d+2)\left(\frac{\wh q}{r}\right)^{d+1}
+r B_0(0) \left(\frac{\wh q}{r}\right)^{d+2},$$
which equals $r B_{d+2}(\wh q/r)$.
The same argument for the coefficient of $\psi^{d+1}$ allows
us to rewrite \eqref{eq:thething2}
as
$$r\frac{B_{d+2}(\wh q/r)\psi^{d+1}+B_{d+2}(q/r)\wh \psi^{d+1}}{\psi+\wh \psi};$$
by
$B_n(x)=(-1)^nB_n(1-x)$ and $-(-1)^{d+1} z^{d+1}+w^{d+1}=(z+w)\sum_{i+j=d} (-z)^iw^j$, we
finally get the term \eqref{eq:lastterm} as desired.
\end{proof}

\begin{rem}[calculation of Chern classes]
Let $\ch_d$ denote the term of degree $d$ of the Chern character.
The $k$th
Chern class is a polynomial in the variables $\ch_d$, see \cite[I/(2.14)']{McD}
for an explicit formula.
In view of the calculation of Chern classes, in Corollary \ref{cor:pushdown},
we derive a formula for $\ch_d$,
allowing us to express all the products of type $\ch_{d_1}\dots \ch_{d_m}$
in terms of known products between tautological classes in
the rational cohomology ring of the moduli stack of $r$-stable
curves $\MMMbar_{g,n}(r)$.
Recall that the coarse space of $\MMMbar_{g,n}(r)$
coincides with that of the Deligne--Mumford compactification $\MMMbar_{g,n}$.
Therefore
the rational cohomology rings are isomorphic and the intersections of
kappa classes, psi classes and boundary classes coincide and are known.
\end{rem}
In order to prepare the statement of Corollary \ref{cor:pushdown},
 we need to introduce a set of
notation analogous to Remark \ref{notn:charclsing}.
\begin{defn}
Let us denote by $V$ the singular locus of the coarse universal
$r$-stable curve and let us write $V'$ for the usual double cover of $V$.
A node $p\in V$ is \emph{nonseparating}
if the normalization  at $p$ of the fibre
$C\ni p$ is connected. We say that $p$ has \emph{type} ${\rm irr}$.

Let $p\in S$ be a separating node in $V$.
Let $p'\in V'$ be a point lifting $p=\ol\varepsilon (p')$.
The point $p'$ is of \emph{type} $(l, I)$, where
$l\in \{0, \dots, g\}$ and $I\subset [n]:=\{1,\dots ,n\}$
if the normalization at $p$ of the fibre
$C\ni p$ is the union of two disjoint
curves $C_1$ and $C_2$ containing
the first and the second branches attached to $p'\in V'$,
having genera $l$ and $g-l$, and
marked at the points $\{x_i\mid i\in I\}$ and
$\{x_i\mid i\in [n]\setminus I\}$,
respectively.

We observe that the type is a locally constant parameter. Therefore we have
$$V'= V'_{\rm irr}\sqcup \textstyle \bigsqcup _{(l,I)}  V'_{(l,I)}.$$
and the natural morphisms
$$i_{\rm irr}\colon V'_{\rm irr}\to \MMMbar_{g,n}(r)\quad \quad \text{and}
\quad \quad
i_{(l,I)}\colon V'_{(l,I)}\to \MMMbar_{g,n}(r).$$
\end{defn}

\begin{rem}
We apply the above definition to the singular locus $Z$ of the universal
$r$-stable curve on $\MMMbar_{g,n}^r$ and we get
\begin{equation}\label{eq:decompZ'1}
Z'= Z'_{\rm irr}\sqcup{\textstyle \bigsqcup _{(l,I)} Z'_{(l,I)}}.
\end{equation}
For any separating node $p$, the normalization at $p$ of the coarse space of
the fibre $C\ni p$
is $C_1\sqcup C_2$, where $C_1$ contains the first branch.
Note that $p'\in Z'$ lifting $p$ belongs to $Z'_q$ according to Notation \ref{notn:charclsing}.(c)
if any only if we have
$$\deg_{C_1}\left((\omega_C^{\log})^{\otimes s}(-\textsum_{m_i}[x_i])\right)
-r\deg _{C_1} (\ol S)=r-q,$$
where $\ol S$ is the line bundle induced by the universal root via pushforward to the
coarse space.
It easily follows that a point $p'\in Z'$ of type $(l,I)$ belongs to $Z_q$
for $q\in \{0, \dots, r-1\}$ satisfying
\begin{equation}\label{eq:qtype}
q+\deg_{C_1}\left((\omega_C^{\log})^{\otimes s}(-\textsum_{m_i}[x_i])\right)\in r\ZZ,
\end{equation}
or equivalently
\begin{equation}\label{eq:qtypeexplicit}
q+ 2ls-s-\textsum_{i\in I}(m_i -s) \in r\ZZ.\end{equation}
In this way we refine the
decomposition of $Z'$ given in  \eqref{eq:decompZ'1}, and we get
the disjoint union
$$Z'={ \textstyle \bigsqcup _{q=0}^{r-1} }Z'_{({\rm irr},q)}\sqcup{\textstyle \bigsqcup _{(l,I)} Z'_{(l,I)}}.$$
We also define the natural morphisms
$$j_{(\irr,q)}\colon Z'_{(\irr,q)}\to \MMMbar_{g,n}(r)\quad \quad \text{and}
\quad \quad
j_{(l,I)}\colon Z'_{(l,I)}\to \MMMbar_{g,n}(r).$$
By definition $Z'_{(l,I)}$ fits in the fibred diagram
\begin{equation}\label{eq:commlI}\xymatrix@R=0.2cm{
Z'_{(l,I)}\ar[rr]^{j_{(l,I)}}\ar[dd]_{p_{(l,I)}}& &\MMMbar_{g,n}^r\ar[dd]^p \\
&\square&\\
V'_{(l,I)}\ar[rr]_{i_{(l,I)}}  &    & \MMMbar_{g,n}(r),}
\end{equation}
whereas, for $j_{(\irr,q)}$, we have the fibre diagram
\begin{equation}\label{eq:commq}\xymatrix@R=0.2cm{
{\textstyle \bigsqcup_{q=0}^{r-1}}Z'_{(\irr,q)}\ar[rr]^{\bigsqcup_q j_{(\irr,q)}}
\ar[dd]_{\bigsqcup _q p_{(\irr,q)}}& &\MMMbar_{g,n}^r\ar[dd]^p \\
&\square&\\
V'_{\irr}\ar[rr]_{i_{(\irr)}}  &    & \MMMbar_{g,n}(r).}\end{equation}
As stated in Theorem \ref{thm:etale}, we have
\begin{equation}\label{eq:plI}\deg(p)=\deg (p_{(l,I)})=r^{2g-1},\end{equation}
because each fibre contains $r^{2g}$ geometric points with stabilizer $\pmmu_r$.
On the other hand, the degree of each morphism $p_{(\irr,q)}$ satisfies
\begin{equation}\label{eq:pqirr}\deg(p_{(\irr,q)})=r^{2g-2}.\end{equation}
This means that the morphisms $p_{(\irr,0)}, \dots, p_{(\irr, r-1)}$
have equal degree.
We show this claim by studying the fibre. Fixing a point $x'\in V_{\irr}'$ means
choosing an
$n$-pointed  $r$-stable curve $C$
and a nonseparating node $x\in C$, which we regard as $B\pmmu_r\into C$.
Furthermore, it means choosing a branch of the node, whose cotangent line can be
regarded as a generator of $\Pic(B\pmmu_r)$. The geometric points of
the fibre of $p_{(\irr,0)}\sqcup \dots \sqcup p_{(\irr, r-1)}$
over $x'\in S'_{\irr}$ are all the
$r^{2g}$ distinct $r$th roots of $\Kcal$
on $C$.
Each point can be assigned to a
nonempty substacks $H_0, \dots, H_{r-1}$
according to
the restriction homomorphism $\Pic(C)\to \Pic(B\pmmu_r)$.
In this way $H_0, \dots, H_{r-1}$
realize the disjoint union of the fibres of $p_{(\irr,0)}$,
\dots, and $p_{(\irr, r-1)}$.
We claim that the stacks $H_i$ are all isomorphic, and therefore their degree
equals $(1/r)\deg p$.
Indeed, we can choose an $r$th root $T_i$
in each substack $H_i$ and, via  $S\mapsto S\otimes T_i^\vee$,
identify $H_i$
to the stack of $r$-torsion line bundles which are trivially linearized  at
the node $x\in C$.
\end{rem}

We denote by $\psi_i$ and $\kappa^{}_d$
the classes defined as in Notation \ref{notn:taut} in the cohomology of $\MMMbar_{g,n}(r)$.
Similarly, we define the classes $\psi$ and
$\wh \psi$ in $V'$. Clearly their pullback via $p$, $p_{(l,I)}$, and $p_{(\irr,q)}$
yields the analogue classes in
$\MMMbar_{g,n}^{r}$. Theorem \ref{thm:GRRcalc} can be written as
follows.
\begin{cor}\label{cor:pushdown}
Let $\SSS$ be the universal $r$th root of
$(\omega^{\log})^{\otimes s}(-\textsum_{i=1}^n m_i[x_i])$
on the universal $r$-stable curve.
We have
\begin{multline}
\left[\ch(R^\bullet \pi_* \SSS)\right]_d=\\
p^*\Bigg(\frac{B_{d+1}}{(d+1)!}\left(s/r\right) \kappa^{}_d -
\sum_{i=1}^n \frac{B_{d+1}}{(d+1)!})\left(m_i/r\right)\psi_i^{d}+
\frac{1}{2}\sum_{\substack {0\le l\le g\\ I\subset [n]}}
r\frac{B_{d+1}(q(l,I)/r)}{(d+1)!}i_{(l,I)*}
\left({\gamma_{d-1}}\right)\Bigg)\\
+\frac{1}{2}\sum_{q=0}^{r-1}r\frac{B_{d+1}(q/r)}{(d+1)!}j_{(\irr,q)*}
\left({\gamma_{d-1}}\right),
\end{multline}
where the cycles $\gamma_{d}$
in $A^{d}(V'_{(l,I)})$ and  $A^{d}(Z'_{(\irr,q)})$
are defined as in Theorem \ref{thm:GRRcalc}
and $q(l,I)$ is the index fitting in \eqref{eq:qtype}. \qed
\end{cor}

\subsection{An example}
With an example
we illustrate how the above formula allows us
to calculate intersection numbers. When the locus $Z'_{\irr}$
of curves with nonseparating nodes  is empty, the
calculation can take place in the
standard  moduli space of
stable curves.  Otherwise, the above
corollary and the degree evaluations \eqref{eq:pqirr} allow us
to carry out the calculation quite easily as
the following example shows.
\begin{exa}[$1$-pointed genus-$1$ curves]\label{exa:11}
In this example we take $s=g=n=1$, $m_1=r+1$.
This means that we look at the stack $\MMMbar_{1,1}^r$ classifying $r$th roots
of $\omega^{\log}(-(r+1)[x_1])$ on $r$-stable genus-$1$ $1$-pointed curves.
In this case the direct image $R^\bullet\pi_* \sta S$ of the universal $r$th root in the derived
category is represented by $-L$ where
$L$ is the line bundle $R^1\pi_* \sta S$; this happens because
each fibre is irreducible and the restriction of
$\sta S$ on every fibre of the universal curve $\pi\colon \sta C\to \MMMbar_{1,1}^r$
has degree $-1$.
As we see in Section \ref{sect:applWspin}, this example
matches Witten's predictions of \cite{Wi} (see Example \ref{exa:inv11}).
We show
\begin{equation}\label{eq:grr11}
\deg {\rm c}_1(R^\bullet \pi_* \sta S)=\deg{\rm c}_1(-L)=-\frac{r-1}{24}.\end{equation}
\end{exa}
First, we use the formula of Corollary \ref{cor:pushdown}
for $d=1$ and get
\begin{equation}
\ch_1= \frac{1}{2}p^*
(B_2\left({1}/{r}\right)\kappa^{}_1-
{B_2\left(1+1/{r}\right)}\psi_1)
+\frac{1}{4}\textsum_{q=0}^{r-1} rB_2(q/r)j_{(\irr, q)*} (\gamma_{0}).
\end{equation}
Since $\pi_*(c_1(\omega_\pi)^2)$
vanishes, we have $\kappa_1=\psi_1$; we write
$$\ch_1= \frac{1}{2}(B_2\left({1}/{r}\right)-{B_2(1+1/r)})p^*\psi_1
+\frac{1}{4}\textsum_{q=0}^{r-1} rB_2(q/r)j_{(\irr, q)*} (\gamma_{0}).
$$
Then the relation $B_n(x)-B_n(x+1)=-nx^{n-1}$ implies
$$\ch_1= \frac{1}{2}(-2(1/r))p^*\psi_1
+\frac{1}{4}\textsum_{q=0}^{r-1} r((q/r)^2-(q/r)+1/6)j_{(\irr, q)*} (\gamma_{0}),
$$
where we used the formula for the second Bernoulli polynomial.
We simplify and we evaluate $\gamma_0=1$
\begin{equation}\label{eq:cross}
\ch_1= -p^*\psi_1/r
+\frac{r}{12}({\textstyle \bigsqcup_q j_{(\irr,q)}})_*(1/2)-
\sum_{q=0}^{r-1} \frac{q(r-q)}{4r}j_{(\irr, q)*}p_{(\irr, q)}^*(1).
\end{equation}
We evaluate the degree via pushforward to $\MMMbar_{1,1}(r)$ (recall the
diagram \eqref{eq:commq})
$$p_* \ch_1= -\psi_1+\frac{r}{12}
i_{(\irr)*}({\textstyle \bigsqcup_q p_{(\irr, q)}})_*(1/2)-
\sum_{q=0}^{r-1}
\frac{q(r-q)}{4r}i_{(\irr)*}p_{(\irr, q)*}p_{(\irr, q)}^*(1).$$
By \eqref{eq:pqirr},
the degree of $p_{(\irr, q)}$ is $1$ (one could also
directly check that each fibre consists of
$r$ copies of $B\pmmu_r$).
We have
$$p_* \ch_1= -\psi_1+\frac{r^2}{12} i_{(\irr)*}(1/2)-
\sum_{q=0}^{r-1} \frac{q(r-q)}{4r}i_{(\irr)*}(1).$$

For any positive integer $r$ we have
\begin{equation}\label{eq:compound}\sum_{q=1}^{r-1}q(r-q)=\frac{(r-1)r(r+1)}{6}.
\end{equation}
Indeed, assume that the relation holds for $r-1$ and write
$\textsum_{q=1}^{r-2}q(r-q-1)=(r-2)(r-1)r/6.$
Adding $\textsum_{q=1}^{r-1}q=(r-1)r/2$ to both sides yields
$\textsum_{q=1}^{r-1}q(r-q)=(r-2)(r-1)r/6 + (r-1)r/2$ and,
therefore, the claim. Using the formula \eqref{eq:compound},
we get
$$p_* \ch_1=
-\psi+\frac{r}{12} i_{(\irr)*}(1/2)-
\frac{r^2-1}{12}i_{(\irr)*}(1/2).$$
Finally, since $V'$ is a double cover of the boundary locus,
we regard $i_{(\irr)*}(1/2)$ as the boundary locus of
$\MMMbar_{1,1}(r)$: the fundamental class of the substack
of singular $r$-stable $1$-pointed genus-$1$ curves. Note that
such a substack contains only one $r$-stable curve up to isomorphism:
the $1$-pointed $r$-stable curve $(X,x)$
over the nodal cubic $(\coa{X}, \coa{x})$
marked at a smooth point.
Therefore, this substack is $B\Aut(X,x)$, where
$\Aut(X, x)$ fits in the natural exact sequence
$$1\to \Aut(({X}, {x}),(\coa{X}, \coa{x}))\to
\Aut({X}, {x})\to \Aut (\coa{X}, \coa{x}) \to 1,$$
where $\Aut(({X}, {x}),(\coa{X}, \coa{x}))$ is the group of automorphisms
fixing the coarse space  and is isomorphic to $\pmmu_r$ by
\cite[Prop.~7.1.1]{ACV}. Clearly
$\Aut (\coa{X}, \coa{x})$ is generated by the
hyperelliptic involution and has order two. Putting everything together we get
$$\deg i_{(\irr)*}(1/2)=1/2r.$$
We evaluate $\deg \psi=1/24$ (this holds for $\MMMbar_{1,1}(r)$, because
the morphism to $\MMMbar_{1,1}$ has degree $1$ and $\psi$ is a pullback
from $\MMMbar_{1,1}$).
We finally reach Witten's intersection number (see Section \ref{sect:applWspin})
$$\deg ( -\psi)+\frac{r}{24} -\frac{r^2-1}{24r}= \frac{-r +r^2 -r^2 +1}{24r}=\frac{1-r}{24r}.$$

\subsection{A different approach via Toen's GRR formula}\label{sect:staGRR}
As mentioned in the introduction, by a result of Toen \cite{To}, the
Grothendieck Riemann--Roch calculation can be carried out
entirely in the context of stacks, rather than via direct image
followed by resolution of singularities.
We illustrate this in detail for the case treated in the previous example.
In Remark \ref{rem:comp}, we briefly illustrate how the calculation can be
carried out in general.

\begin{exa}
Consider the moduli stack $\MMMbar_{1,1}^r$ as in the
previous example. Let $\wt \pi\colon \wt C\to X$ be
the $r$-stable curve of genus $1$ on
a scheme $X$ \'etale on the moduli stack.
Let $\Delta$ be the divisor specifying the universal marking.
The universal $r$th root is a line bundle defined on $\wt C$.
Let $\wt Z$ be the singular locus of $\wt C\to X$; note that
$\wt Z$ has dimension $0$ and it is a  union of copies of $B\pmmu_r$.
Let $\wt Z'$ be the double cover of $\wt Z$ classifying the
orders of the branches; as in Remark \ref{rem:singloc}, we have
\begin{equation}\label{eq:Ssplitq}
\wt Z'= {\textstyle \bigsqcup _{q=0}^{r-1}} \wt{Z}'_q\end{equation}
The morphism $\wt \pi$ is a morphism between a regular
Deligne--Mumford stack  and a regular scheme; therefore
Toen's theorem applies. In order to state it, we need to introduce some
constructions deriving from the
initial setting.
\begin{itemize}
\item[$\bullet$]
First, note that all stabilizers in $\wt C$
are isomorphic to
$\pmmu_r$. Throughout the example,
$\xi$ is a primitive $r$th root of unity.
\item[$\bullet$]
We consider the inertia stack $I$ of $\wt C$. By definition, it is
the fibre product of $\wt C$ with itself over $\wt C\times \wt C$.
In this case, we have the following explicit description:
we take
a copy $\wt C{({0})}$ of
$\wt C$ and $r-1$ copies $\wt Z({1}), \dots, \wt Z({r-1})$ of $\wt Z$, and
we have
\begin{equation}
\label{eq:inertia}
I= \wt C{({0})}\sqcup {\textstyle \bigsqcup _{i=1}^{r-1}\wt Z({i}}),
\end{equation}
together with the natural morphism $\eta\colon I\to \wt C$.
Naturally, as soon as we pass to the
double cover $\wt Z'({i})$
each copy decomposes as a disjoint union of
substacks $\wt Z_q'(i)$ and we have the morphisms
\begin{equation}\label{eq:jqi}
j_{q,i}
\colon \wt Z_q'(i)\into \wt Z'(i)\to I,
\end{equation}
analogue to the morphisms introduced in Notation
\ref{notn:charclsing}.(g).
The above morphisms are defined for $0\le q<r$ and $0< i< r$; however, by
taking a further copy
$\wt Z(0)$ and $\wt Z'(0)$
of $\wt Z$ and $\wt Z'$ we extend the above
definition \eqref{eq:jqi} to the indices $i=0$ and $0\le q<r$ so that
$j_{q,0}$ maps to $I$ factoring through $C(0)$.
\item[$\bullet$]
In the $K$-theory ring of $I_{\wt C}$ tensored with $\CC$, we consider the
homomorphism
\begin{equation}\label{eq:rhodiag}
\rho\colon K^0(I)\otimes {\CC}\to K^0(I)\otimes {\CC}\end{equation}
defined as follows: if a bundle $W$ is decomposed into
a direct sum of eigenbundles $W^{(i)}$ with eigenvalue $\xi^i$
\begin{equation}\label{eq:rho}
W={\textstyle \bigoplus_{i=0}^{r-1} W^{(i)}};
\quad \quad \text{then, we have} \quad \quad
\rho[W]=\sum_{i=0}^{r-1}\xi^i [W^{(i)}].\end{equation}
\end{itemize}

Toen's GRR formula reads
\begin{equation}\label{eq:toen}\ch(R^\bullet\wt \pi_* \wt S) = \wt\pi_*\eta_* \left(\ch (
\rho(\eta^*\wt S)) \ch\left(\rho(\la_{-1}\mathcal N_{\eta}^\vee )\right)^{-1}
\td^\vee( \Omega_{\wt \pi\circ \eta})\right),\end{equation}
where $\la_{-1}$ is
$\sum_i(-1)^i\bigwedge^i$ and we use
the fact that $\ch(\rho(\la_{-1}\mathcal N_{\eta}^\vee ))$ is
invertible
\cite[Lem.~4.6]{To}.

On the connected component labelled with $0$ in the inertia stack,
the formula between brackets on the right hand side yields
$$\frac{K\exp(K/r)}{1-\exp(K)} \exp(-\Delta) +
\frac{1}{12}\textsum_{q=0}^{r-1}(j_{q,0})_*(1/2),$$
where $K$ is the first Chern
class of the relative canonical line bundle.
The degree $1$ part yields the usual Riemann--Roch calculation
of the degree of $R^\bullet \pi_*\sta S$.
The contribution in degree $1$ on $X$ comes from the
classes of degree $2$.
After pushforward to $X$ we get
$$-\frac{\psi_1}{r}
+\frac{1}{12r}({\textstyle \bigsqcup_q j_{q}})_*(1/2).$$

On the connected components
$\wt Z'_q{(i)}$ for $i=1,\dots, r-1$ the
formula \eqref{eq:toen}
reads
\begin{equation}\label{eq:toenoni}
\frac{\xi^{qi}}{(1-\xi^i)(1-\xi^{-i})}.
\end{equation}
where $\xi^{qi}$ comes from the character
induced by the universal $r$th root on $\wt Z'_{q}(i)$ and
from the definition \eqref{eq:rho} of $\rho$,
whereas $(1-\xi^i)(1-\xi^{-i})$ is
obtained after splitting $\mathcal N^\vee$ into two lines
with $\pmmu_r$-linearization $\xi\mapsto \xi$ and
$\xi\mapsto \xi^{-1}$, applying the definition of $\rho$,
and using
the fact that $\la_{-1}$ is multiplicative.

A routine calculation using the formula $1/(1-\xi)=
\sum_{d=i}^{r-1} (d/r) \xi^i$
yields
\begin{equation}\label{eq:cyclolittle}
\sum_{i=1}^{r-1}\frac{\xi^{qi}}{(1-\xi^i)(1-\xi^{-i})}=\frac{r^2-1}{12}
-\frac {q(r-q)}{2}.
\end{equation}
Therefore,
by pushing forward all the
classes of dimension $0$ mentioned above to $X$, we get
$$-\frac{\psi_1}{r}
+\frac{1}{12r}({\textstyle \bigsqcup_q j_{q}})_*(1/2)
+\frac{r^2-1}{12r}({\textstyle \bigsqcup_q j_{q}})_*(1/2)
-\sum _{q=1}^{r-1} \frac{q(r-q)}{2r}({\textstyle \bigsqcup_q j_{q}})_*(1/2).$$
Now, we can continue the calculation just as in the previous example,
after identifying the above equation with \eqref{eq:cross}
\end{exa}

\begin{rem}\label{rem:comp}
We illustrate how this calculation can be carried out in general.
As in the example above the calculation
consists of a first part involving the component
$\wt C(0)$ in the
inertia stack of $\wt C$ and
a second part involving the remaining components
$\wt C(1),\dots, \wt C(r-1)$.
We use the same notation as in the example.

On $\wt C(0)$, the evaluation of the product $\ch (
\rho(\eta^*\wt S)) \ch\left(\rho(\la_{-1}\mathcal N_{\eta}^\vee )\right)^{-1}
\td^\vee( \Omega_{\wt \pi\circ\eta})$ appearing in Toen's formula
\eqref{eq:toen} can be obtained as follows.
First, on $\wt C(0)$, the universal $r$th root
is a line bundle $\wt S$ and an $r$th root
of $(\omega^{\log})^{\otimes s}(-\sum_i m_i \Delta_i)$. In this way, we have
$$\ch(\wt S)= \exp (sK/r){\textstyle \prod_{i=1}^n }\exp ((s-m_i)\Delta_i/r).$$
Note also that the normal sheaf $\mathcal N_\eta$ vanishes on $\wt C(0)$.
The formula expressing the Todd character is the same as in
Mumford's calculation. Therefore, we can proceed as in Step 2 of the proof of the main theorem
and we get the summands involving the kappa classes and the psi classes.

The evaluation of $ \ch\left(\rho(\la_{-1}\mathcal N_{\eta}^\vee )\right)^{-1}
\td^\vee( \Omega_{\wt \pi\circ\eta})$
on the components $\wt C(1),\dots, \wt C(r-1)$ can be found in
\cite[\S7.2.6]{Ts} for any twisted curve. In the case of
$r$-stable curves, the expression involves
cohomology classes with coefficients in the $r$th cyclotomic field
rather than $\QQ$.
Finally, $\ch (
\rho(\eta^*\wt S))$ can be evaluated as follows.
Note first that the restriction of the universal $r$th root $\wt S$
on $\wt C(1), \dots, \wt C(r-1)$ is an $r$th root of $\Ocal$ ($\omega$
and its twists at the markings are
trivial on the singular locus). This allows us to show that, on $\wt C(i)$,
the term $\ch(\rho(\eta^*\wt S))$
can be regarded as the evaluation on $\xi^i$ of the
character associated to $\wt S$.
After these evaluations we can check
that the total contribution in the cohomology ring of the base scheme
yields the last summand of the formula of Theorem \ref{thm:GRRcalc}.
To this effect, as in the example above,
we need to express the sum of the terms with coefficients in the representations
in terms of rational coefficients. In the example,
we used the equation \eqref{eq:cyclolittle}, whereas the general
result can be obtained by applying the formula
$$\sum_{i=1}^{r-1} \frac{\xi^{qi}}{\xi^{i}e^x-1}=\frac{re^{(r-q)x}}{e^{rx}-1}-\frac{1}{e^x-1}.$$
\end{rem}

\section{Applications and motivations}\label{sect:appl}
In the recent years the
interest in the enumerative geometry of $r$th roots has been revived by
the Gromov--Witten theory of orbifolds and Witten's
$r$-spin variant of Gromov--Witten theory. We briefly illustrate
some cases where Theorem \ref{thm:GRRcalc}
can be applied: the geometry of roots of $\Ocal$ are related to
the crepant resolution conjecture in the context of
Gromov--Witten theory of orbifolds, whereas the roots of
$\omega$ are related to Witten's $r$-spin theory.

\subsection{Enumerative geometry of the orbifold $[\CC^2/\pmmu_r]$}\label{sect:applCRC}
As mentioned in the introduction, the crepant resolution conjecture for $\CC^2/\pmmu_r$
provides the most transparent example of an
application of Theorem \ref{thm:GRRcalc}
in Gromov--Witten theory.
We discuss this prototype example following \cite{BGP}; we only work with curves
of genus $0$. In view of Coates and Ruan's recent conjectures \cite{CR} mentioned in the introduction,
one can easily adapt the algorithm illustrated here
to higher genus.

As in \cite{BGP} we
adopt a choice of a primitive $r$th root
of unity $\xi$ all throughout the example.
To begin with, let us recall a known fact on twisted curves (see \cite{Ca}, \cite{Ol}).
\begin{rem}[twisted curves]\label{rem:cad}
We recall that the notion of twisted curves is more general
than the notion of $r$-stable curve: it allows stabilizers
of any finite order on the nodes
and on the markings, see \cite{AV}.
Adding stabilizers on the markings is easy:
for any positive multiindex $l_1, \dots, l_n$,
there is the following equivalence of categories.
On the one hand, we consider the category of
$n$-pointed twisted curves
with representable (scheme-theoretic) smooth locus.
On the other hand of the equivalence, we
have the category of $n$-pointed twisted curves whose $i$th marking
has stabilizer of order $l_i$, see \cite[Thm.~4.1]{Ca}.

We show how this equivalence can be used.
Let $C\to X$ be the morphism with representable smooth locus and
write as $C(l_1,\dots, l_n)\to X$
the twisted curve corresponding to $C$
in the equivalence.
Then, $C(l_1, \dots, l_n)$ carries a
tautological $l_i$th root $M_i$ of
$\Ocal(D_i)$ where $D_i$ is the divisor specifying
the $i$th marking. We also have a natural  morphism
of stacks $C(l_1, \dots, l_n)\to C$.
Any line bundle $L$ on $C(l_1, \dots, l_n)$ can be described uniquely as
a pullback from $C$ tensored with $M_i^{\otimes m_i}$ for $0\le m_i<l_i$,
see \cite[Cor.~2.12]{Ca}.
In this way, any line bundle on $C(l_1, \dots, l_n)$
induces local indices $m_i\in \{0, \dots, l_i-1\}$ at each node.
We say that $L$ is \emph{faithful} if $m_i$ is prime to $l_i$.
We finally recall the following fact:
the pushforward via $C(l_1, \dots, l_n)\to C$ identifies
faithful $r$-torsion line bundles on the twisted curve $C(l_1, \dots, l_n)$
with indices  $m_i$ on the $i$th point with $r$th roots
of $\Ocal(-\sum_i r m_i/l_i D_i)$ on  $C$.
\end{rem}

\paragraph{The statement.}
We start from the $\pmmu_r$-action
on $\CC^2$
$$\xi\cdot(x,y)\mapsto (\xi x, \xi^{-1}y)$$
which commutes with the action of the torus $T=\mathbb G_m\times \mathbb G_m$
on $\CC^2$.
The statement relates
the equivariant Gromov--Witten potentials of the resolution $Y$ of the
scheme $\CC^2/\pmmu_r$
and that of the stack $\sta X=[\CC^2/\pmmu_r]$.
In the most recent formulation \cite[Conj.~10.2]{CR}
the conjectural statement is that the two data can entirely
identified
after a change of variables.

The definition of the equivariant Gromov--Witten potential of $Y$ is classical,
see \cite[\S2]{BGP}.
Here, we illustrate how Theorem \ref{thm:GRRcalc} applies to the calculation of
the Gromov--Witten potential of $\sta X$. We set $g=0$ for simplicity.
We recall the definition of the intersection numbers appearing in the
equivariant genus-$0$ Gromov--Witten potential of $\sta X$.
The $T$-action on $\CC^2$ induces a $T$-action on $\sta X$.
Consider the moduli stack of stable maps to $\sta X$.
The inertia stack $I_{\sta X}$
has $r$ components corresponding to the elements of $\pmmu_r$, and
each component is contractible; so, the cohomology of the inertia stacks has a canonical basis
$\{D_0, D_1, D_2, \dots, D_{r-1}\}$, where $D_0$ is
the class of the point at the origin in the special component of
the inertia stack identified to $\sta X$.
The equivariant Gromov--Witten potential of $\sta X$
is just a power series whose coefficients are classes in
the ring of $T$-equivariant cohomology of the point
$${\rm GW}(n_0, \dots, n_{r-1})\in H_{T}(\pt)=\QQ[t_1,t_2].$$
This is the precise definition of ${\rm GW}(n_0, \dots, n_{r-1})$
in $T$-equivariant cohomology
\begin{equation}\label{eq:GWpot}
\int _{[{\MMMbar}_{0, n_0+\dots +n_{r-1}}(\sta X,0)]^{\rm vir}}\prod_{i=1}^{n_0} {\rm ev_i}^* (D_0)
\prod_{i=n_0+1}^{n_0+n_1} {\rm ev_i}^* (D_1)
\dots
\prod_{i=n_0+\dots +n_{r-2}+1}^{n_0+\dots+n_{r-1}} {\rm ev_i}^* (D_{r-1}),
\end{equation}
where ${\MMMbar}_{0, n_0+\dots +n_{r-1}}(\sta X,0)$ is the stack of
stable maps to $\sta X$ and $\rm ev_i$ denotes the $i$th evaluation map (see \cite{AV}
and \cite{AGV}).

In \eqref{eq:GWred}, we write the formula \eqref{eq:GWpot} more explicitly
in order to make clear that Theorem \ref{thm:GRRcalc} can
be efficiently applied in order to calculate the invariants.
We focus on the case $n_0=0$ and $\sum_i n_i>0$, which is
the difficult part of the calculation.
We proceed in three steps. In Step 1, we
describe more explicitly the space ${\MMMbar}_{0, n_0+\dots +n_{r-1}}(\sta X,0)$.
In Step 2, we point out that the integral \eqref{eq:GWpot} only involves
a special connected component of  ${\MMMbar}_{0, n_0+\dots +n_{r-1}}(\sta X,0)$.
In Step 3, we make explicit what we mean by integration against the virtual class.

\begin{enumerate}
\item \text{Step 1. The moduli stack.}\\ Since $n_0=0$, stable maps of degree 0
factor through
$B\pmmu_r\subset \sta X$; therefore we can integrate over
$\MMMbar_{0, n_0+\dots+n_{r-1}}(B\pmmu_r)$.
We recall that $\MMMbar_{0, n_0+\dots+n_{r-1}}(B\pmmu_r)$ is equipped with a
universal curve and a universal morphism
$$\pi\colon \sta C\to \MMMbar_{0, n_0+\dots+n_{r-1}}(B\pmmu_r) \quad \quad \text{and}\quad \quad
f\colon \sta C\to B\pmmu_r,$$
where $f$ is representable.
The above universal object is in fact a twisted curve of genus $0$
equipped with an $r$-torsion line bundle.
The condition of \emph{representability} imposed on
$f$ can be rephrased as saying that the line bundle is faithful.
Therefore, for each marking,
there exist two \emph{coprime} indices $(l,m)$
$$\text{$l$ \ \ \ \ and \ \ \ \ $0\le m< l$,}$$
such that the local picture is
the line bundle $L_U$ on $[U/\pmmu_l]$, with $U=\Spec \CC[z]$, $\pmmu_{l}$-action
given by $z\mapsto \xi_{l}z$, and $L_U$ equal to $U\times \CC$ with
$\pmmu_l$-linearization $\chi\colon \xi_{l}\mapsto \xi_{l}^{-m}$.
Note that in this way we have defined
the multiindices $(m_1, \dots, m_n)$ and $(l_1, \dots, l_n)$, which  are
locally constant on $\MMMbar_{0, n_0+\dots+n_{r-1}}(B\pmmu_r)$.

\item\text{Step 2. The class
$\prod_{i=1}^{n_0} {\rm ev_i}^* (D_0)
\dots\prod_{i=n_0+\dots +n_{r-2}+1}^{n_0+\dots+n_{r-1}} {\rm ev_i}^* (D_{r-1})$.}\\
It is easy to see that such a class is supported on the connected component
$$Z(n_0,n_1, \dots, n_{r-1})\subset \MMMbar_{0, n_0+n_1+\dots+n_{r-1}}(B\pmmu_r)$$ of stable maps whose
local indices $(l_i,m_i)$ at the $i$th point equal $(r/\hcf\{r,j\},j/\hcf\{r,j\})$,
where
$j$ is the integer in $\{0,1, \dots, r-1\}$ satisfying
$n_0+\dots + n_{j-1}+1 \le i\le n_0+\dots+ n_j$ (set $n_{-1}=0$).
Using Remark \ref{rem:cad}, the stack
$Z(n_0,n_1, \dots, n_{r-1})$
can be regarded as the stack of faithful $r$th roots of
\begin{equation}
\label{eq:Kcrc}\textstyle {\Kcal=\Ocal \left(-\sum_{i=n_0+1}^{n_0+n_1} [x_i] - \sum_{i=n_0+n_1+1}^{n_0+n_1+n_2} 2[x_i]
-\dots- \sum_{i=n_0+n_1+\dots+n_{r-2}+1}^{n_0+n_1+\dots+n_{r-1}} (r-1)[x_i] \right)}
\end{equation}
on the universal twisted curves
with trivial stabilizers at the markings.

\item\text{Step 3. The virtual class.} \\
Finally, for
dimension reasons we impose $n_0+n_1+\dots+ n_{r-1}>3$ (which
amounts to requiring $n_1+\dots+ n_{r-1}>3$ since $n_0$ vanishes).
Furthermore, for degree reasons we notice that
the connected component $Z(0,n_1, \dots, n_{r-1})$
is empty unless we have
$$\textsum_i in_i\in r\ZZ$$ (this condition specializes to
$n_1\equiv n_2 \mod 3$ in the case $r=3$, see \cite[A.2.(3)]{BGP}).
In the remaining cases the equivariant Gromov--Witten invariant is possibly
nonzero and the virtual class can be evaluated as
the Euler $T$-equivariant class of $R^1 \pi_*f^*
(L_\xi\oplus L_{\ol\xi})$
$$e_T(R^1 \pi_*f^*
(L_\xi\oplus L_{\ol\xi})),$$ where
$L_\xi$ and $L_{\ol \xi}$ are the line bundles on $B\pmmu_r$
associated to the $\pmmu_r$-representations $\xi\mapsto \xi$ and $\xi\mapsto \xi^{-1}$.
Notice that $f^*L_\xi$ and $f^*L_{\ol \xi}$
on $Z(0,n_1, \dots, n_{r-1})$ can be regarded as
the universal $r$-torsion line bundle $\sta T$ of $\Ocal$ on $\sta C$ and its dual
line bundle
$\sta T^\vee$, respectively.
For $n=\textsum _in_i>3$ and
$\textsum_i i n_i\in r\ZZ$, we get
\begin{multline}\label{eq:GWred}
{\rm GW}(0, n_1,\dots, n_{r-1})= \deg e_T(R^1 \pi_*f^*
(L_\xi\oplus L_{\ol\xi}))= (t_1+t_2)\deg {\rm c}_{n-3}(R^1\pi_* (\sta T \oplus \sta T^\vee)),\end{multline}
where the last equality is a calculation in $T$-equivariant cohomology
based on Mumford's proof of $c(\mathbb H\oplus \mathbb H^\vee)=0$, see \cite[Lem.~3.1]{BGP}.
\end{enumerate}

We finally draw the connection to the universal
$r$th root on $r$-stable curves and to Theorem \ref{thm:GRRcalc}.
The calculations \eqref{eq:GWred} can be carried out on
the stack $\MMMbar_{0,n}^{r}$ of $r$th roots on
$n$-pointed $r$-stable curves
of the line bundle $\Kcal$ of \eqref{eq:Kcrc}.
It should be noted that the stack $\MMMbar_{0,n}^{r}$ of $r$th roots
on $r$-stable curves differs from $Z(n_0, n_1,\dots, n_{r-1})$, because in the
definition of $\MMMbar_{0,n}^{r}$ we require
that all stabilizers at the nodes equal $r$. However,
there is a natural morphism
invertible on the interior, whose ramification indices
on the boundary
can be explicitly handled, see \cite[Thm.~4.1.6]{Ch_mod}.

On $\MMMbar_{0,n}^{r}$,
the pullback of $R^1\pi_* (\sta T \oplus \sta T^\vee)$
yields the
direct image via the universal $r$-stable curve $\pi$ of
the sum of $\sta S $ and $\sta S^\vee(-\sum_i [x_i])$.
Note that for both line bundles the sections of the
restriction to each fibre of the $r$-stable curve $\pi$ are necessarily zero
because their $r$th power vanishes on each component (this is an easy consequence of
Remark \ref{rem:cad} or
of Lemma \ref{lem:Dq} from the scheme-theoretic point of view).
In this way, Theorem \ref{thm:GRRcalc} provides the Chern character
of $R^1\pi_* (\sta T \oplus \sta T^\vee)$, whereas
the standard formulae
expressing Chern classes in terms of Chern characters \cite[I/(2.14)']{McD}
allow us to derive the Gromov--Witten invariants.
In general this calculation, although effective, is long to carry out
in practice for high values of $n$ and $g$.
For $g=0$, the formula of Corollary \ref{cor:pushdown} can be further
simplified because $Z_{\rm irr}$ is empty. In general, one can
easily implement
this algorithm by adapting Faber's computer programme for the calculation of the
intersection numbers of tautological classes.

\begin{exa}
We calculate the genus-$0$ Gromov--Witten invariant corresponding to
$r=5$ and $$(n_0,n_1,n_2,n_3,n_4)=(0,0,3,0,1).$$ This amounts to consider
the moduli space of curves of genus $0$ with $4$ markings (the sum $\sum_in_i$) equipped with
a $5$th root of $\Ocal(-2[x_1]-2[x_2]-2[x_3]-4[x_4])$. The points of this moduli functor
can be identified with the points of
the projective line $\MMMbar_{0,4}\cong \PP^1$, and the points of the boundary
can be regarded as $3$ points in $\PP^1$. However,
the automorphisms are nontrivial.
Each point representing smooth curves has  stabilizer $\pmmu_5$ (the
multiplication by a root of unity along the fibres of the root); furthermore,
the boundary points have a further automorphisms acting about the
 node as $(x,y)\mapsto (x,\xi y)$ (see \cite[Prop.~7.1.1]{ACV}).
In this way, this moduli functor can be regarded as a
stack whose coarse space is the projective line with $25$ automorphisms
on $3$ points forming the boundary locus $\Delta$ and
$5$ automorphisms on all remaining points.

In order to calculate $\text{GW}(0,0,3,0,1)$, we integrate the first Chern class of
the vector bundle $R^1\pi_*\mathcal S\oplus R^1\pi_*\mathcal S^\vee(-\sum_i[x_i])$,
where $\mathcal S$ is the universal root of
$\Ocal(-2[x_1]-2[x_2]-2[x_3]-4[x_4])$. This amounts to evaluating
the sum of the degrees of
$\ch_1(R^1\pi_*\mathcal S))$ and $\ch_1(R^1\pi_*\mathcal S^\vee(-\sum_i[x_i]))$. As noted above
$\mathcal S^\vee(-\sum_i[x_i])$ can be regarded as a universal root of
$\Ocal(-3[x_1]-3[x_2]-3[x_3]-[x_4])$. Then, Theorem \ref{thm:GRRcalc} yields an
identity between the terms of degree $1$ of the Chern character of
$R^1\pi_*\mathcal S$ and $R^1\pi_*\mathcal S^\vee(-\sum_i[x_i])$ (this fact
can be generalized and simplifies all calculations).
Applying Corollary \ref{cor:pushdown}, we finally get
\begin{align*}
\text{GW}(0,0,3,0,1)=&2\ch_1(R^1\pi_*\mathcal T)\\
         =&2\left(\frac{B_2(0)}{2}\kappa_1-\frac{B_2(2/5)}{2}\sum_{i=1}^3\psi_i-\frac{B_2(4/5)}{2}\psi_4
         +\frac{5B_2(1/5)}{2}\Delta\right)\\
         =&2\left(\frac{1}{12}\kappa_1+\frac{11}{300}\sum_{i=1}^3\psi_i-\frac{1}{300}\psi_4
         +\frac{5}{300}\Delta\right)\\
         =&\frac{2}{5}\left(\frac{1}{12}+\frac{33}{300}
         -\frac{1}{300}+\frac{3}{300}\right)\\
         =&\frac{2}{5}\left(\frac{1}{12}+\frac{35}{300}\right)=\frac{2}{25}\\
\end{align*}
Indeed, one should notice that the Bernoulli coefficient $q(0,I)$ for all subsets of
the set of markings is equal to $1$ or $4$. Since $B_2(q/r)$ does not vary if we replace
$q$ by $5-q$, all points of $\Delta$
appear with the same coefficient.
\end{exa}
\begin{exa}
Similarly, we get
\begin{align*}
\text{GW}(0,3,1,0,0)=&2\ch_1(R^1\pi_*\mathcal T)\\
         =&2\left(\frac{B_2(0)}{2}\kappa_1-\frac{B_2(1/5)}{2}\sum_{i=1}^3\psi_i-\frac{B_2(2/5)}{2}\psi_4
         -\frac{5B_2(2/5)}{2}\Delta\right)\\
         =&2\left(\frac{1}{12}\kappa_1-\frac{1}{300}\sum_{i=1}^3\psi_i+\frac{11}{300}\psi_4
         -\frac{55}{300}\Delta\right)\\
         =&\frac{2}{5}\left(\frac{1}{12}-\frac{3}{300}
         +\frac{11}{300}-\frac{33}{300}\right)\\
         =&\frac{2}{5}\left(\frac{1}{12}-\frac{25}{300}\right)=0\\
\end{align*}
We point out that for $r=3$
the
invariants GW$(0,n_1, n_2)$
only depend on the sum
of the entries $n_i$, \cite[Prop.~A.1]{BGP}.
As we see here,
this is not the case in general.
\end{exa}

\subsection{Enumerative geometry of $r$-spin structures}\label{sect:applWspin}
By definition $r$-spin structures are the objects of the moduli stacks
$\MMMbar_{g,n}^r$ when we set $s=1$, $2g-2+n>0$, and $2g-2-n-\sum_im_i\in r\ZZ$.
This section illustrates the connection between the
theory of $r$-spin curves and Theorem \ref{thm:GRRcalc}.
The reader can refer to \cite{CZ} for a result
showing the role of Theorem \ref{thm:GRRcalc} in the more general context of
Gromov--Witten theory of $r$-spin structures.

We proceed as follows: we define the $r$-spin potential
$F^{\text {$r$-spin}}$,
we provide a slight reformulation of the definition
of $F^{\text {$r$-spin}}$ (Proposition \ref{pro}), we illustrate with some examples
(Examples \ref{exa:inv11}--\ref{exa:torsion2}) the
 relation with Theorem \ref{thm:GRRcalc}, and we briefly discuss
further conjectural applications of Theorem \ref{thm:GRRcalc} to
 the enumerative geometry of
Hurwitz numbers.

\paragraph{The potential $F^{\text {$r$-spin}}$.}
In \cite{Wi}
Witten introduces the potential $F^{\text {$r$-spin}}$, which
involves the construction of
$\MMMbar_{g,n}^r$, the definition of the
Witten top Chern class ${\rm c}_{\rm W}$, and
the intersection of ${\rm c}_{\rm W}$ with the psi classes.

The moduli stacks $\MMMbar_{g,n}^r$ is the stack of $r$th
roots of $$\mathcal K(\vec m)=\omega^{\log}(-\textsum _i m_i [x_i]).$$
This definition makes sense for nonnegative indices $g,n$ satisfying $2g-2+n>0$ and
any multiindex $\vec{m}=(m_1, \dots, m_n)$
satisfying
$2g-2+n - \sum _i m_i\in r\ZZ$.

The Witten top Chern class ${\rm c}_{\rm W}(\vec{m})$ is a cohomology class of
degree $-\chi(\sta S(\vec{m}))=h^1(\sta S(\vec{m}))-h^0(\sta S(\vec{m}))$ in
$\MMMbar_{g,n}^{r}$ (see \cite{PV} or \cite{Ch_K} for
two different but equivalent definitions). Its definition makes sense under the above conditions
on $g, n$ and $\vec m$ and a further condition $m_i>0$ for any $i$.

We define
$${\rm W}_{g,n}(\vec{a}, \vec{m})=
\int _ {\MMMbar_{g,n}^{r}} \psi_1^{a_1}\dots \psi_n^{a_n}{\rm c}_{\rm W}(\vec{m})\in\QQ.$$
The definition of $F^{\text {$r$-spin}}$ only involves multiindices $\vec m$ satisfying $0<m_i\le r$
for any $i$.

\begin{defn}[the potential $F^{\text {$r$-spin}}$] The potential $F^{\text {$r$-spin}}$
is the power series
$$F^{\text {$r$-spin}}=\sum_{\substack{g,n\ge 0\\2g-2+n>0}} \ \sum_{\lvec{a}, \lvec{m}}
\frac{{\rm W}_{g,n}(\lvec{m},\lvec{a})}{{n}!}
\prod_{i=1}^{n} m_i(m_i+r)\dots(m_i+r(a_i-1)) (m_i+ra_i)\  t_{m_i+ra_i},
$$
where $\vec{m}$ and $\vec{a}$  are multiindices with $n$ coordinates
and the summation is taken over
the pairs $(\lvec{a}, \lvec{m})$ satisfying
$0< \lvec{m}\le r$ and $2g-2+n-\sum_i m_i\in r\ZZ$.
\end{defn}

\begin{rem}\label{rem:donotneed}
The integral of a cycle
of positive dimension vanishes. We do not need to restrict to
the multiindices $\vec {a}$ and $\vec{m}$ for which we have
$\deg({\rm c}_{\rm W}(\vec m))+\sum _i a_i=\dim \MMMbar_{g,n}^r.$
\end{rem}

\begin{rem}
The power series $F^{\text {$r$-spin}}$ can be regarded as the Gromov--Witten
$r$-spin potential of a point and is generalized in \cite{JKVmaps} to the case where
the target is a K\"ahler manifold.
\end{rem}

\begin{thm}[\cite{FSZ} (Witten's $r$-spin conjecture \cite{Wi})]
The power series $\exp(F^{\text {$r$-spin}})$ is a tau function for the Gelfand--Diki\u\i\ hierarchy
$r$-KdV. \qed
\end{thm}
\begin{rem}
The statement above determines $F^{\text {$r$-spin}}$ and
allows new calculations in the enumerative geometry of curves
equipped with $r$-spin structures.
Indeed, by Theorem \ref{thm:GRRcalc}, in \cite{CZ}
we provide a tool for deriving the
twisted Gromov--Witten $r$-spin potentials from $F^{\text {$r$-spin}}$.
\end{rem}

\paragraph{Reformulation of the definition of $F^{\text {$r$-spin}}$.}
Using a result of Polishchuk and Vaintrob (see Lemma \ref{lem:desc})
conjectured by Jarvis, Kimura, and Vaintrob \cite{JKV2}, we
reformulate the above definition of $F^{\text {$r$-spin}}$
via a change of variables.
(This change of variables is also suggested in \cite[Rem.~2.3]{JKV2}.)

Recall Remark \ref{rem:donotneed}: we can either sum over all indices
$\vec m$ and $\vec a$ for which $d$ divides the degree $2g-2+n-\sum_i m_i\in r$,
or we can select those  multiindices $\vec m$ and $\vec a$ for which
the dimension of $\MMM_{g,n}^r$ is equal to $\deg({\rm c}_{\rm W}(\vec m))+\sum _i a_i$.
It turns out that this second choice has several advantages. Let us make this condition explicit.

It is easy to see that $\vec m$ and $\vec a$ satisfy
$\deg({\rm c}_{\rm W}(\vec m))+\sum _i a_i=\dim \MMMbar_{g,n}^r$
if and only if \begin{equation}\label{eq:ccc}\textsum_i(ra_i+m_i)=(2g-2+n)(r+1).\end{equation}

\begin{rem}
An immediate advantage of imposing \eqref{eq:ccc} directly in the
definition of $F^{\text {$r$-spin}}$ is that
we do not need to impose the condition
$2g-2+n-\sum_i m_i\in r\ZZ$: it follows as an immediate consequence.
\end{rem}

\begin{rem}
A less obvious advantage of imposing \eqref{eq:ccc} is that
each number ${\rm W}(\vec a, \vec m)$ can be computed without involving
psi classes. We see how in two steps.

First, we change parameters.
The equation \eqref{eq:ccc} suggest that the two multiindices $\vec m$
and $\vec a$ can be incorporated
in a single multiindex $$\vec{k}=r\vec a +\vec m.$$
Instead of starting from two multiindices $\vec m$ and $\vec a$ satisfying \eqref{eq:ccc},
we can start from a single multiindex $\vec k$ satisfying
\begin{equation}\label{eq:kccc}
\textsum_i k_i =(2g-2+n)(r+1).
\end{equation}
The datum of a positive multiindex $\vec k$ and that of
two multindices  $\vec m\in \{1,\dots, r\}^n$ and $\vec a\in (\ZZ_{\ge 0})^n$
are interchangeable.

Second, we take $\vec{k}=r\vec a +\vec m$ and
we calculate ${\rm W}_{g,n}(\vec m, \vec a)$
without psi classes.
On the one hand, recall that by definition ${\rm W}_{g,n}(\vec m, \vec a)$ is the cup product
$\psi_1^{a_1}\dots \psi_n^{a_n}{\rm c}_{\rm W}(\vec{m})$.
On the other hand compute ${\rm W}_{g,n}(\vec k, \vec 0)$;
explicitly, this means:
take the moduli stack of $r$th roots
of $\Kcal(\vec k)$
and apply the construction of ${\rm c}_{\rm W}$ to the
universal $r$th root $\sta S(\vec{k})$ (as pointed out above this
construction makes sense as soon as the multiindex is positive and not only
inside $\{1,\dots, r\}$).
Since $\vec k$ satisfies \eqref{eq:kccc}
we get a cycle ${\rm c}_{\rm W}(\vec k)$ of dimension $0$
and therefore, by integrating, a number.
The following lemma shows how ${\rm W}_{g,n}(\vec m, \vec a)$ can be computed
from ${\rm W}_{g,n}(\vec k, \vec 0)$.
\end{rem}
\begin{lem}[{\cite[Prop.~5.1]{PV}}]\label{lem:desc}
The classes ${\rm c}_{\rm W}(\vec{m})$ and
${\rm c}_{\rm W}(\vec{k})={\rm c}_{\rm W}(r\vec{a}+\vec{m})$ satisfy
\begin{equation*}\label{eq:desc}{\rm c}_{\rm W}(\vec{k})=\left(\prod_{i=1}^{n}
\frac{ m_i\dots(m_i+(a_i-1)r)}
{r^{a_i}}\right)
\psi_i^{a_1} \dots \psi_n^{a_n}{\rm c}_{\rm W}(\lvec{m}).\qed \end{equation*}
\end{lem}
\begin{cor}\label{cor:desc}
For $\vec{m}$ and $\vec{a}$  satisfying
\eqref{eq:ccc} and $0<m_i\le r$, we have
$${\rm W}_{g,n}(\vec{m},\vec{a})={\left(\prod_{i=1}^{n}
\frac{r^{a_i} }
{m_i\dots(m_i+(a_i-1)r)}\right)} {{\rm W}_{g,n}(0,\vec{k})},$$
where $\vec k$ equals $r\vec a+\vec m$.\qed
\end{cor}
We get the following simplified definition
of the potential $F^{\text {$r$-spin}}$.
\begin{pro}\label{pro}
We assume as usual $r\ge 2$ and  $g,n$ satisfying
$(2g-2+n)>0$.
We have
$$F^{\text {$r$-spin}}=\sum_{\substack {g,n\ge 0\\ 2g-2+n >0}} \ \sum_{\substack{k_i>0 \\ {\sum_i k_i=(2g-2+n)(r+1)}}}
\frac{\deg {\rm c}_{\rm W}(\vec{k})}{n!}
\prod_{i=1}^{n} m_ir^{\lfloor m_i/r\rfloor}t_{m_i},
$$
where ${\rm c}_{\rm W}(\vec{k})$ is
Witten's top Chern  associated to the universal $r$th root of
$\omega^{\log}(-\sum_i k_i[x_i])$.\qed
\end{pro}

\paragraph{The Witten top Chern class and the direct image {$R^\bullet\pi_*{\sta S(\vec{k})}$}.}
We recall from \cite{PV} and \cite{Ch_K} that the
class ${\rm c}_{\rm W}(\vec{k})$
is determined by the direct image
$R^\bullet \pi_* \sta S(\vec k)$ and by
a natural homomorphism in the derived category induced by $\sta S(\vec{k})^{\otimes r}\cong \Kcal(\vec{k})$ and Serre duality:
\begin{equation}\label{eq:form}\Sym^{r}(R^\bullet \pi_* \sta S)\to \Ocal [-1].\end{equation}
Such a class is a generalization of the top Chern class:
it coincides with $\ctop (R^1\pi_*\sta S)$ when $\sta S$ has no sections along each fibre
of the universal $r$-stable curve (note that this implies $\pi_*\sta S=0$ and $R^1\pi_*\sta S\in \Vect$).
In general ${\rm c}_{\rm W}$
is not the top Chern class of a vector bundle, but it
has the same cohomological degree as $-\chi(\sta S(\vec{k}))$

\begin{rem}We point out that
for $\lvec{k}$ satisfying
$\textsum_{i=1}^n k_i=(r+1)(2g-2+n)$
we have $-\chi(\SSS)=\dim(\MMMbar_{g,n}^r)$. Furthermore
$\sum _i k_i =(r+1)(2g-2+n)$ implies
that the direct image $\pi_*L$ vanishes, because the relative degree
of the universal $r$th root of $\Kcal(\vec{k})$ is negative:
$$\deg_{\pi}\SSS(\vec{k})=\frac{2g-2+n -\textsum_{i=1}^n m_i}{r}= {-(2g-2+n)}<0.$$
On the open substack of
$r$th roots defined on \emph{irreducible} $r$-stable curves
the direct image $R^1\pi_*L$ is a vector bundle
of rank $-\chi(\SSS)=\dim(\MMMbar_{g,n}^r)$,
because the $r$-spin structure has no section
on each fibre. Away from the locus classifying irreducible
curves, the direct image is not a vector bundle in general.
We show this with some examples.
\end{rem}

Assuming $\sum_ik_i=(r+1)(2g-2+n)$, we provide three examples.
In the first case, $R^1\pi_*\sta S$ is locally free and Witten's
top Chern class is its top Chern class.
In the second and third example,
$R^1\pi_*\sta S$ has
torsion (even if the condition $\sum_ik_i=(r+1)(2g-2+n)$
guarantees  that $\pi_*\sta S$ vanishes the sheaf $R^1\pi_*\sta S$ is not locally free).
The second example shows that this happens in codimension $1$ inside
any moduli space $\MMMbar_{g,n}^r$ for $g>0$ and $(g,n)\ne (1,1)$.
The third example has the advantage of being very elementary and concrete,
stack-theoretic curves do not appear.

\begin{exa}\label{exa:inv11}
Since every $1$-pointed $1$-stable curve is irreducible,
we deduce that $\deg {\rm c}_{\rm W}(r+1)$ can be regarded as the first Chern class of
the line bundle $R^1\pi_*\sta S$ on the stack
$\MMMbar_{1,1}^r$ of $r$th roots of $\omega_\pi (-r[x_1])$. In this way, using
the calculation carried out in Example \ref{exa:11} and
the above Lemma \ref{lem:desc},
we get a new verification of the prediction of Witten \cite{Wi}
$$\int _{\MMMbar_{1,1}^r} \psi_1 {\rm c}_{\rm W}(1)={\rm W}_{1,1}(1, 1)=
r{\rm W}_{1,1}(0, r+1)= \frac{r-1}{24},$$
see \cite[4.4]{JKV}, \cite[5.3/3]{PV}, and \cite[6.4]{Ch_K}
for different verifications of the same invariant.
\end{exa}

\begin{exa}\label{exa:torsion1}
We exhibit an $r$-spin structure on a curve
with $n$ markings and multiindex $\vec{k}$ satisfying $\sum_ik_i=(r+1)(2g-2+n)$.
Consider a curve $C$ of genus $g$ with one node and two
smooth irreducible components of genus $1$ and $g-1$.
Assume that the markings belong to the component of genus $g-1$.
Note that, for degree reasons,
all $r$-spin structures are defined on the curve having stabilizer of order $r$
on the node.
Furthermore, since all markings are on the component of genus $g-1$,
the sections of any $r$-spin structure
constantly vanish on the component of genus $g-1$.

We are interested in the group of
sections $H^0$ of these $r$-spin structures. We want to show that
even if $\sum_ik_i=(r+1)(2g-2+n)$ is satisfied, it may well happen that $H^0$
does not vanish.
Using Lemma \ref{lem:Dq}, one can regard the sections of any $r$-spin structure
on $C$ as sections of
a line bundle $S'$ on a nodal curve $C'$
obtained from $C$ by normalization and by gluing a projective line to the
two preimages of the node.
There exists an $r$-spin structure for
which the corresponding line bundle $S'$ on $C'$
is simply obtained by gluing
three line bundles defined on each component and  verifying the following conditions:
the line bundle on the component of genus $1$ is trivial,
the line bundle on the rational component is isomorphic to $\Ocal(1)$,
the line on the remaining component has only one global section $s=0$.
Now it is immediate to see that for such a
line bundle $H^0$ is one-dimensional.
\end{exa}
\begin{exa}\label{exa:torsion2}
Set $r=2$,
consider a curve of genus $g>1$ with an even number $n$ of
markings, and choose a multiindex $k_1, \dots, k_n$
such that $\sum _i k_i= 3(2g-2+n)$.
If the curve is smooth, the $2$-spin structures are square roots of a line bundle of
degree $-2(2g-2+n)$, therefore for all $2$-spin structures we have $H^0=0$.
We let the curve degenerate and acquire $2$ nonseparating nodes: we now look at a curve
$C$ which is the union of a projective line $P$
and a curve $E$ of genus $g-1$ meeting at two points.
We choose the nodal curve in such a way that
the markings belong to the irreducible component of genus $0$ (in this way the nodal curve is stable). We consider
a square root of $\omega^{\log}_C(-\sum k_i[x_i])$. Its degree on $P$ is
$-\sum _i (k_i-1)/2=-3(g-1)-n$, whereas
on $E$ the degree is $g-1$. Since the degree on
$P$ is strictly negative,
the sections of a square root $S$ in $\Pic C$, are sections of
the restriction on $E$ vanishing at the points meeting $P$.
By construction, the restriction $S\rest E$ is a theta characteristic of the curve
$E$ of genus $g-1$. We can easily choose from the start a curve
$E$  with a theta characteristic having more than two
linearly independent sections (see Mumford, \cite[\S4]{Mu_theta} for explicit examples).
\end{exa}

The problem of finding an explicit closed formula for
Witten's $r$-spin invariants is still open. Using the
cohomological field theory axioms \cite{JKV}, Shadrin and Zvonkine provide an
algorithm \cite{SZ}. More generally,
it would be interesting to know
how the Witten's top Chern class is related to the tautological classes
(see \cite[Conj.~2.5.3]{FJR} for a conjectural statement) and to
the Chern classes of
$R^\bullet\pi_*\sta S$ (see the discussion in \cite{Po}).
\paragraph{The Chern classes of $R^\bullet \pi_* \sta S$ and Hurwitz numbers.}
Recently, pursuing a different approach to Witten's statement,
Zvonkine formulated
new conjectures on the enumerative
geometry of the $r$-spin structures. He starts
from the Ekedahl--Lando--Shapiro--Vainshtein formula
of {\cite{ELSV}}:
\begin{equation}\label{eq:ELSV}
h_{g;b_1,\dots ,b_n} =
(d + n + 2g - 2)!
\prod_{i=1}^n
\frac{b_i^{b_i}}{b_i!}
\int _{\MMMbar_{g,n}}
\frac{1 - \la_1 + \la_2 - \dots {\pm} \la_g}
{(1 - b_1 \psi_1) \dots (1 - b_n \psi_n)},\end{equation}
where $d=\sum_i b_i$, $\lambda_k$ denotes the $k$th Chern class of the Hodge bundle,
and $h_{g;b_1,\dots ,b_n} $ is the Hurwitz number of simply ramified covers of the sphere
with profile $(b_1, \dots, b_n)$ at infinity.
His conjectural formula involves in a
similar way the Chern classes
${\rm c}_k(R^\bullet\pi_*\sta S)$ (see \cite{ZvLu}).

\noindent
Email: \url {chiodo@ujf-grenoble.fr}\\
Address: Institut Fourier\\
U.M.R. CNRS 5582\\
U.F.R. de Math\'{e}matiques\\
Universit\'{e} de Grenoble 1\\
BP 74, 38402\\
Saint Martin d'H\`{e}res\\
France}

\end{document}